\documentclass[11pt]{article}
\usepackage{epsfig}
\usepackage{amsmath}
\usepackage{amsfonts}
\usepackage{amssymb}

\title{Jakobson's Theorem near saddle-node bifurcations}

\author{Ale Jan Homburg\\
KdV Institute for Mathematics \\
University of Amsterdam \\
Plantage Muidergracht 24 \\
1018 TV Amsterdam \\
The Netherlands \\
{\small \ttfamily alejan@science.uva.nl}
\and
Todd Young \\
Department of Mathematics \\
Morton Hall \\
Ohio University \\
Athens, OH 45701 \\
U.S.A.\\
{\small \ttfamily young@math.ohiou.edu}}

\newcounter{bean}

\newtheorem{main}[bean]{Theorem}

\newtheorem{theorem}{Theorem}[section]

\newtheorem{lem}[theorem]{Lemma}
\newtheorem{prop}[theorem]{Proposition}

\newtheorem{defn}[theorem]{Definition}

\newtheorem{rem}[theorem]{Remark}

\newcommand{\comp}{\lower3pt\hbox{${}^{\circ}$}}

\newfont{\bb}{msbm10  scaled\magstep1}
\newfont{\bbs}{msbm8  scaled\magstep0}
\newfont{\euf}{eufm10 scaled\magstep1}
\newfont{\efm}{eufm8  scaled\magstep0}

\newcommand{\Res}[2]{{#1}\raisebox{-.4ex}{$\left|\,_{#2}\right.$}}
\newcommand{\qed}{\hspace*{\fill}$\square$}

\newcommand{\bN}{\mbox{\bb N}}
\newcommand{\bR}{\mbox{\bb R}}
\newcommand{\bZ}{\mbox{\bb Z}}

\newcommand{\Ac}{{\mathcal A}}

\newcommand{\Dc}{{\mathcal D}}

\newcommand{\Gc}{{\mathcal G}}

\newcommand{\Lc}{{\mathcal L}}
\newcommand{\Oc}{{\mathcal O}}
\newcommand{\Pc}{{\mathcal P}}
\newcommand{\Qc}{{\mathcal Q}}

\newcommand{\Uc}{{\mathcal U}}

\newcommand{\bS}{{\bf S}}

\newcommand{\Mb}{{\bf M}}
\newcommand{\Nb}{{\bf N}}

\newcommand{\Rb}{{\bf R}}

\newcommand{\ga}{\gamma}

\newcommand{\ep}{\epsilon}
\newcommand{\ka}{\kappa}
\newcommand{\la}{\lambda}

\newcommand{\Om}{\Omega}

\newcommand{\De}{{\Delta}}

\newcommand{\ft}{{\tilde{f}}}
\newcommand{\Et}{{\tilde{E}}}

\newcommand{\taub}{{\bar{\tau}}}
\newcommand{\ct}{\tilde{c}}
\newcommand{\BR}{\mathrm{BR}}

\begin{document}

\maketitle

\begin{abstract}
We discuss one parameter families of unimodal maps,
with negative Schwarzian derivative, unfolding a
saddle-node bifurcation. 
It was previously shown that for a parameter set of positive
Lebesgue density at the bifurcation, the maps possess
attracting periodic orbits of high period.
We show that there is also a parameter set of positive 
density at the bifurcation, for which  
the maps exhibit absolutely continuous 
invariant measures which are supported on the largest possible
interval. We prove that these measures converge weakly to
an atomic measure supported on the orbit of the saddle-node point.
Using these measures we analyze the intermittent time series that result
from the destruction of the periodic attractor 
in the saddle-node bifurcation and prove 
asymptotic formulae for the frequency with which 
orbits visit the region previously occupied by the periodic attractor.
\end{abstract}

\section{Introduction}\label{sec_intro}

\subsection{Background}

This article is a companion article for \cite{homyou00}.
In that article we discussed intermittent dynamics 
associated with boundary crisis (homoclinic) bifurcations in families of
unimodal maps. In the present work we treat 
saddle node bifurcations from the same perspective.

By Jakobson's celebrated work \cite{jak81}, the logistic family
$x \mapsto \mu x (1-x)$, 
$x \in [0,1]$,  admits absolutely continuous invariant measures (a.c.i.m.'s)
for $\mu$ from a set of positive measure.
In fact, $\mu = 4$ is a (Lebesgue or full) density point of this set.
A different argument for this result was given by
Benedicks and Carleson \cite{bencar85}, \cite{bencar91}.
Their reasoning was generalized to unfoldings
$\{ f_\ga\}$ 
of unimodal Misiurewicz maps $f_0$, with
eventually periodic critical point $c$ (and possessing 
negative Schwarzian derivative). 
It was shown that the bifurcation value $\ga =0$ is a density point
of the set of parameter values for which $f_\ga$ admits
an absolutely continuous invariant measure 
\cite{melstr93}, \cite{thitreyou94}.

We generalize these results to unfoldings of saddle-node bifurcations
in families $\{ f_\ga \}$  of unimodal maps
with negative Schwarzian derivative.
We establish that a saddle-node bifurcation value occurs as a point
of positive density of the parameter set for which there are
absolutely continuous invariant measures. 
Under the assumption that $f_0$ is not more than once renormalizable,
we construct parameters for which $f_\ga$ possesses a.c.i.m.'s supported on the maximal interval
$[f_\ga^2 (c) , f_\ga(c)]$, see Theorem~\ref{existence} below.
In contrast to the Misiurewicz bifurcation values, 
the saddle-node bifurcation value
is not a full density point of this parameter set.
This is because it is known that the parameter set for which
the map has a periodic attractor has positive density at
any saddle-node.

Following the construction of a.c.i.m.'s, we continue with
a detailed discussion of the intermittency that
occurs due to the saddle-node bifurcation.
That saddle-node bifurcations can give rise to intermittency
is known since \cite{pomman80}, who called intermittency 
associated with a saddle-node bifurcation type I intermittency.
Pomeau and Manneville studied type I intermittency 
in connection with the Lorenz model.
In the model, 
simplifying (hyperbolicity) assumptions on the dynamics
outside a neighborhood of the saddle-node periodic orbit are made.
In perhaps the most basic example of type I intermittency, 
in families of unimodal maps, such simplifications are not justified.
This is due to the presence of a critical point.
Our discussion of absolutely continuous 
invariant measures allows us to give a rigorous
treatment of intermittent time series, where we explain and prove 
quantitative aspects earlier discussed numerically 
in \cite{hirhubsca82}, see Theorems~\ref{intermittency} and 
\ref{saddlenode} below.

Diaz et.\ al.\ \cite{diarocvia96} studied the unfoldings of
saddle-node bifurcations in higher dimensional diffeomorphisms
and their results imply in the present context that there
exists a subset $\Om_D$ which has positive density at $0$, 
such that for each $\ga \in \Om_D$, 
$f_\ga$ has an absolutely
continuous invariant measure. The measures
produced there are supported on small periodic domains on 
which the map renormalizes to a H\'{e}non-like family.
In the present context this corresponds to parameter values inside
periodic windows, for which $f_\ga$ is renormalizable.
However, the invariant measures
produced in Jakobson's work are 
supported on the maximal possible interval, 
$[f^2_\ga(c),f_\ga(c)]$ and the invariant measures that we 
construct are also
supported on $[f^2_\ga(c),f_\ga(c)]$. 

While editing the final draft of this paper, we learned of the 
results of Maria Jo\~ao Costa \cite{cos01}, corresponding to part of 
her 1998 thesis, on related work in a similar context.
She focused on the sink-horseshoe bifurcation in which a sink 
and a horseshoe collapse, see \cite{zee82,cos98}, 
and studied families $\{g_\ga\}$
of unimodal interval maps to describe bifurcations.
The class of families of interval maps studied in \cite{cos01} 
differs from ours; it consists of
unfoldings of unimodal maps $g_0: [0,1] \to [0,1]$, that 
possess a saddle node fixed point at $p \in (0,1)$, so that the critical
point $c > p$ satisfies $g_0^2 (c) < p$.
For such families Costa derived the analogue of Theorem~\ref{existence}
below, stating that absolutely continuous invariant measures supported on the
maximal interval occur with positive density at the bifurcation point.
The set-up of Costa's proof corresponds to ours in that she 
also combines
arguments originating from Benedicks and Carleson  \cite{bencar85,bencar91}
with the introduction of induced maps, as discussed below.

\subsection{Assumptions and statement of main results}

Let $\{f_\gamma\}$ be a family of unimodal maps of the interval $[0,1]$,
with critical point at $c$.
Suppose that each $f_\ga$ is at least $C^3$ smooth and that
$f_\ga(x)$, $Df_\ga(x)$, and $D^2 f_\ga(x)$, are $C^1$ w.r.t.\ $\ga$. 
Suppose that each $f_\ga$ has negative Schwarzian 
derivative (see \cite{melstr93}) and that $D^2f_\ga(c) <0$. 
Further, suppose that $f_\ga (1) = f_\ga(0) = 0$ and that 
the fixed point at 0 is hyperbolic repelling.
We say that $\{f_\gamma\}$ {\em unfolds a (quadratic) saddle-node} if,
\begin{itemize}
\item There is a $q$-periodic point $a$, with $Df^q_0 (a) = 1$, and,
\item $D^2 f^q_0 (a) \frac{\partial}{\partial \gamma} f_\gamma^q (a)> 0$ 
at $\gamma = 0$.
\end{itemize}
For the sake of clarity, we will assume that
$$
    \frac{\partial}{\partial \ga} f_\ga^q(a)|_{\ga = 0} > 0 
      \quad \text{and} \quad  D^2f^q_0(a) > 0.
$$
With this convention, for $\ga > 0$ the saddle-node point 
disappears and complicated dynamics may occur.

For a set $A$ of parameter values, let $m(A)$ be its Lebesgue 
measure. A unimodal map $f$ is called {\em renormalizable} if there 
is a proper subinterval $I \subset [0,1]$ containing the critical 
point $c$, so that $f^n(I) \subset I$ for some $n>0$ and 
$f^i(I) \cap I = \emptyset$ for $0<i<n$. A renormalizable map is 
called {\em once renormalizable} if the above property
defines $n$ uniquely. 

\begin{main}\label{existence}
Let $\{f_\ga\}$ be as above, unfolding a saddle-node bifurcation
at $\ga = 0$. Assume that $f_0$ is once renormalizable.
Consider the set $\Gamma$ such that for each $\ga \in \Gamma$,
the map $f_\ga$ has an absolutely continuous invariant measure
$\nu _\ga$ whose support is the 
maximal interval $[f^2_\ga(c),f_\ga(c)]$.
The set $\Gamma$ has positive density at $\ga = 0$:
\begin{eqnarray*}
\liminf_{\ga^\ast\searrow 0} 
\frac{m( \Gamma \cap (0,\ga^\ast))}{ \ga^\ast} & > & 0.
\end{eqnarray*}
It does not have full density:
\begin{eqnarray*}
\limsup_{\ga^\ast\searrow 0} 
\frac{m( \Gamma \cap (0,\ga^\ast))}{ \ga^\ast} & < & 1.
\end{eqnarray*}
\end{main}

\begin{figure}[htp]

\centerline{\hbox{\epsfig{figure=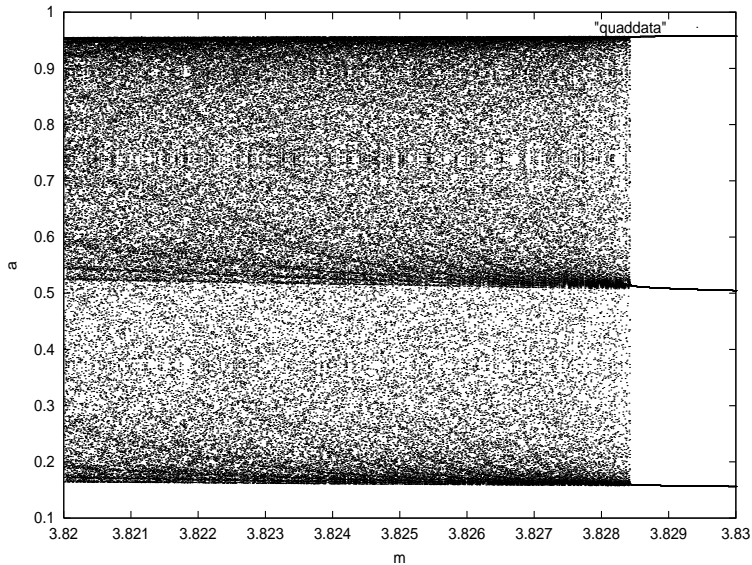,height=3in,width=4in}}}


\caption{Bifurcation diagram for the quadratic family $x \mapsto \mu x (1-x)$.
         For $\mu = \mu_{sn} = 1 + 2\sqrt{2}$ the map undergoes a 
         saddle-node bifurcation of a period 3 orbit. For most observed 
         parameter values, $\mu < \mu_{sn}$, orbits appear to fill the 
         interval $[f^2_\mu(1/2),f_\mu(1/2)]$.}
\label{fig:quad}

\end{figure}

Careful numerical studies of the quadratic family,
$x \mapsto \ga x (1-x)$, predict that  
between the period doubling limit and $\ga = 4$, excluding 
the period three window, less then about 15\% 
of the parameter values correspond to the periodic windows \cite{brian}.
Numerical simulations also suggest that near a saddle-node bifurcation,
there is a large set of parameter values for which the a.c.i.m.'s  
are supported on the maximal possible interval.
This is illustrated in Figure~\ref{fig:quad}. 

Existence of the parameter set $\Gamma$ follows from
a similar approach as used by Benedicks and Carleson 
in their treatment of Jakobson's result.
Below we comment in more detail on the construction of the
parameter set $\Gamma$. We now first consider intermittency that
results from the saddle-node bifurcation and introduce our
second main result.

Intermittent dynamics manifests itself 
by alternating phases
with different characteristics. In one phase, referred to as the 
laminar phase, the dynamics appear to be nearly periodic. 
While in the other phase, the relaminarization phase, 
the orbit makes large, seemingly chaotic
excursions away from the periodic region. These excursions
are called chaotic bursts.
Let $\bar{E}$ be a neighborhood of 
the orbit, $\Oc(a)$ of $a$, not containing a critical point of $f^q_0$.
Let $\chi_{\bar{E}}$ be defined as 
\begin{equation}\label{defn_Phi}
\chi_{\bar{E}}(x,\ga)  =  
\lim_{n\to \infty} \frac{1}{n} \sum_{i=0}^{n-1} 1_{\bar{E}} (f^i_\gamma(x)),
\end{equation}
whenever the limit exists, where $1_{\bar{E}}$ is the usual indicator 
function of the set ${\bar{E}}$.
That is, $\chi_{\bar{E}}(x,\ga)$ is the relative frequency with 
which the orbit ${\cal O}(x)$ visits $\bar{E}$
(for those $x$ for which the limit exists).
The following theorem discusses $\chi_{\bar{E}}(x,\ga)$ for $\gamma$ near 0.
\begin{figure}[ht]

\centerline{\hbox{\epsfig{figure=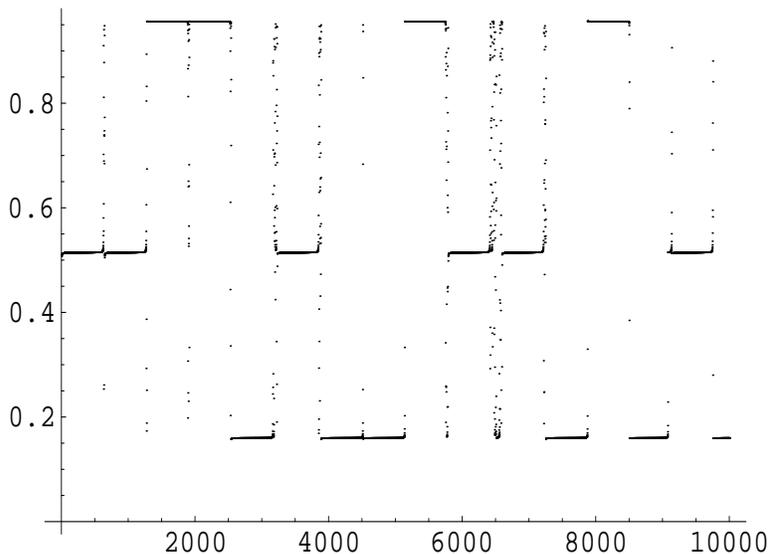,height=3in,width=4in}}}


\caption{Time series for the quadratic map
$x \mapsto \mu x (1-x)$ for $\mu = 3.828$ near the saddle-node
bifurcation of the period three orbit. Both the laminar phase 
(nearly periodic)
and chaotic bursts are clearly seen.}
\label{fig_sntime}

\end{figure}

Denote by $\nu_0$ the atomic measure supported on the orbit
of $a$, given by
\begin{equation}\label{nu0}
     \nu_0 = \frac{1}{k} \sum^{k-1}_{i=0} \delta_{f_0^i(a)}.
\end{equation}
We will denote the usual weak convergence of measures by the symbol 
$\rightharpoonup$.

\begin{main}\label{intermittency}
Let $\{f_\ga\}$ and $\Gamma$ be as in Theorem~\ref{existence}.
There exist sets $\Omega\subset \Gamma$ of parameter values with positive density at $\ga =0$,
so that 
\begin{eqnarray*}
\lim_{\ga \in \Omega, \ga \searrow 0} \nu_\ga 
         & \rightharpoonup & \nu_0.
\end{eqnarray*}
Restricting to $\gamma \in \Omega$, $\chi_{\bar{E}}(x,\gamma)$ 
is a constant, $\chi_{\bar{E}}(\ga)$, almost everywhere on $[0,1]$ and
$\chi_{\bar{E}}(\ga)$ depends continuously on $\gamma$ at $0$.
There exist $K_1,K_2 > 0$ so that 
\begin{equation}\label{pomeau}
K_1 \le \lim_{\gamma\in \Omega, \gamma\searrow 0} 
\frac{1 - \chi_{\bar{E}}(\ga)}{\sqrt{\gamma}}
 \le  K_2.
\end{equation}
\end{main}

\begin{figure}[ht]

\centerline{\hbox{\epsfig{figure=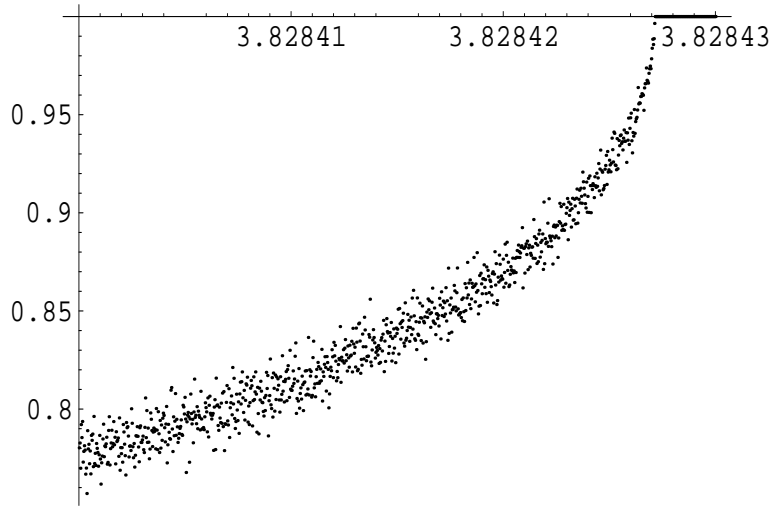,height=3in,width=4in}}}


\caption{Numerical computation of $\chi_{\bar{E}}(\mu)$ as a function of $\mu$ 
for the logistic family 
$f_\mu (x) = \mu x (1-x)$ near the saddle-node bifurcation of a period three orbit
at $\mu_{sn} = 1+2\sqrt{2}$.} 
\label{fig:freq}

\end{figure}
%
%
%
%
Each set $\Omega$ from the above theorem is constructed 
as a union $\cup_{l\ge l_0} \Omega_l$ of sets $\Omega_l$,
where each $\Omega_l$ has, as a Lebesgue density point, 
a parameter value  $\ga_l^*$ for which
$f_{\ga_l^*}^{q l +m } (c)$ hits a hyperbolic repelling periodic point 
$y^*_{\ga_l^*}$,
for some $m$ (the map $f_{\ga_l^*}$ is a Misiurewicz map).
Different sets $\Omega$ result from considering
points $y^*_{\ga_l^*}$ in different periodic orbits.
For the construction of $\Omega_l$, we adapt Benedicks and Carleson's
treatment of Jakobson's result (we follow in fact
the approach taken by Luzzatto \cite{luz00})
and construct parameter values
for which $f_\ga (c)$ has a positive Lyapunov exponent.
The main obstacle in the 
present case is the fact that when orbits fall into
the saddle-node region they remain there for many iterations and the derivative
of $f^q_\ga$ in that region is close to one. We  get around this problem by
defining an induced map which skips over
the saddle-node region. We show that the 
induced map has uniform exponential growth of derivatives and 
that this implies
exponential growth of derivatives (with weaker constants) for the 
original map. 
Careful estimates bound the measure of the sets $\Omega_l$ for 
increasing values of $l$, thus bounding the measure of their union
$\Omega$.

For comparison we include a result from \cite{homyou00},
which builds on results in \cite{afrliuyou96}, \cite{diarocvia96}, 
and \cite{afryou98}, showing that also periodic attracting 
orbits are found for parameters from a set with positive 
density at the bifurcation point. The frequency with which 
the dynamics is in the laminar phase behaves
in a similar way to the frequency in Theorem~\ref{intermittency}.
Note that the fact that the set $\Gamma$ in Theorem~\ref{existence}
does not have full density at $\ga = 0$, follows from the following result.

\begin{main}\label{saddlenode}
Let $\{f_\ga\}$ be as above, unfolding a saddle node bifurcation 
at $\ga = 0$.
There exist sets $\Ac$ of parameter values of positive measure
and positive density at 0, i.e.
\begin{eqnarray}\label{sndensity}
\lim_{\ga \searrow 0}
\frac{m\left( \Ac \cap [0 , \ga)\right) }
{\ga} & > & 0,
\end{eqnarray}
so that for each $\ga \in \Ac$,
$f_\ga$ has an attracting periodic orbit.
Further, \begin{eqnarray*}
\lim_{\ga \in \Ac, \ga \searrow 0} \nu_\ga 
         & \rightharpoonup & \nu_0,
\end{eqnarray*}
where $\nu_\ga$, $\ga \in \Ac$, is the invariant measure 
supported on the periodic orbit.
Restricting to $\gamma \in \Ac$, $\chi_{\bar{E}}(x,\gamma)$ 
is a constant, $\chi_{\bar{E}}(\ga)$, almost everywhere on $[0,1]$ and
$\chi_{\bar{E}}(\ga)$ depends continuously on $\gamma$ at $0$.
There exists $K>0$ so that 
\begin{eqnarray*}
\lim_{\gamma\in \Ac, \gamma\searrow 0} 
\frac{1 - \chi_{\bar{E}}(\ga)}{\sqrt{\gamma}}
& = & K.
\end{eqnarray*}
\end{main}

Different sets $\Ac$ lead to different limit values $K$ in Theorem~\ref{saddlenode}. In fact, the proof of 
Theorem~\ref{saddlenode} makes clear that 
arbitrary large numbers occur as the limit values $K$.
This fact, together with Theorems~\ref{existence}~and~\ref{saddlenode} 
lead us to conjecture that
there is a parameter set $\Lambda$ which has $\ga = 0$ 
as a Lebesgue density point, so that 
$$
 \lim_{\ga\searrow 0,\ga\in\Lambda} \frac{\ln (1 - \chi_{\bar{E}})}{\ln \ga}  = \frac{1}{2}.
$$
It was shown in \cite{homyou00} that such a limit cannot hold without restricting the parameter set.


\section{The saddle node: local embedding flows}
\label{sec_snlocal}

Denote by $E$ a small neighborhood of 
$a$ on which $f^q_0$ is invertible. Let $W^s_{loc}(a)$ and 
$W^u_{loc}(a)$ denote the usual local stable and local unstable sets for
$a$. 

\begin{prop}\label{embedding}
Let $\{f_\ga\}$ be a $C^1$ family of $C^r$, $r \ge 2$, maps unfolding a
saddle-node. Then there
exists a family of $C^r$ flows, $\{\phi^t_\ga\}_{0 \le \ga < \bar{\ga}}$, 
on $E$ such that
$f^q_\ga \equiv \phi^1_\ga$ for each $\ga \ge 0$. Further,
$\phi^t_\ga(\cdot) \rightarrow \phi^t_0(\cdot)$ in the $C^1$ topology
on $E$ and in the $C^r$ topology on compact intervals away from the 
fixed point. The flow $\phi^t_0$ is uniquely determined by $f_0$.
\end{prop}

\noindent {\sc Proof.}
The $C^\infty$ version of this theorem is due to Takens \cite{tak73}. 
The $C^r$ result follows from Part 2 of \cite{yoc95}.
The case $\ga = 0$ follows from Appendix 3 of \cite{yoc95}.
The case $\ga >0$ and the convergences as $\ga \searrow 0$ follow from 
Theorem~IV.2.5 and Lemma~IV.2.7 of the same. 
\qed

\begin{rem}\label{IL}
This result is known as the Takens Embedding Theorem.
A version of it appears in \cite{ilyli99}. They proved that one
may obtain $\phi^t_\ga(x)$ which depends $C^r$ smoothly on both $x$ and
$\ga$, even at the fixed point, if one requires that 
$(x,\ga) \mapsto f_\ga(x)$ be $C^{R(r)}$ smooth, where 
$R(r)$ may be larger than $r$. 
Proposition~\ref{embedding} allows for our weaker hypotheses and its
implications are sufficient for our purposes.
\end{rem}

Choose $d \in W^s_{loc}(a)$ and let
$$
I^s_\ga = [d,f_\ga^q(d)].
$$
Also, choose a point $e \in W^u_{loc}(a)$ so that
$$
I^u_\ga  = [e,f^q_\ga(e)] \subset E.
$$
For the sake of convenience we restrict $E$ to be the
interval
$$
   E = [d,e],
$$
so that $W^s_{loc}(a) = [d,a)$ and $W^u_{loc}(a) = (a,e]$.
We will use the embedding flow on the interval 
$E \cup I^u_\ga = [d,f^q_\ga(e)]$.


Given $\ga \ge 0$ and $x \in I^u_\ga$, define $\tau_\ga^u(x)$ 
to be the unique number for which 
$$
\phi_\ga^{\tau^u_\ga(x)}(e) = x.
$$ 
For $\ga \ge 0$ and $x \in I^s_\ga$, let $\tau_\ga^s(x)$ be defined by 
$$
\phi_\ga^{\tau_\ga^s(x)}(d) = x.
$$
It follows from the smoothness of $\phi^t_\ga(x)$ that for each $\ga \ge 0$,
the functions $\tau_\ga^{s,u}$ are $C^r$ diffeomorphisms from $I^{s,u}_\ga$
to [0,1]. We will use $\tau_\ga^{s,u}$ as coordinates on $I^{s,u}_\ga$.
In the following, we will associate $[0,1]$ with the unit circle $\bS^1$
by making the identification $0 \sim 1$. We also treat $I^{s,u}_\ga$
as circles through the coordinates $\tau^{s,u}_\ga$.

Given $d$ and $e$ as above, let $\{\ga_l\}_{l=l_0}^\infty$
be the sequence, $\bar{\ga} > \ga_{l_0} > \ga_{l_0+1} > \ldots$, defined by:
$$
    f^{ql}_{\ga_l}(d) = e.
$$
For each $l \ge l_0$ let 
$g_l:[0,1] \rightarrow [\ga_{l+1},\ga_l]$,
be the reparameterization map defined by
$$
\phi^{l+\theta}_{g_l(\theta)}(d) = e.
$$
We have that $g_l(0) = \ga_l$ and $g_l(1) = \ga_{l+1}$.
We may invert $g_l(\cdot)$, for each $l$, to obtain maps
$\theta_l:[\ga_{l+1},\ga_l] \rightarrow [0,1]$.

\begin{prop}\label{distortion}
The reparameterization maps $g_l$ are smooth monotone decreasing functions
with uniformly small distortion:
given $\varepsilon > 0$ there is $N \in \bN$ so that for
every $l \ge N$ and every $\theta \subset [0,1]$,
$$
(1-\varepsilon) \le \frac{ Dg_l(\theta)}{|\ga_l - \ga_{l+1}|}
\le (1+\varepsilon).
$$ 
\end{prop}

\noindent {\sc Proof.}
 Diaz et.\ al.\ proved this result under the hypothesis that
$(x,\ga) \mapsto f_\ga(x)$ is $C^{R(r)}$, using Il'yashenko and Li's
embedding result (Remark~\ref{IL}).
A proof of this result under the current hypotheses appears 
in \cite{afryou98} based on \cite{jon90}.
\qed \\

We remark that $l^2 \ga_l$ converges as $l \rightarrow \infty$ (see \cite{MK90}), so that 
$\ga_{l+1}/\ga_l \to 1$ as $l\to \infty$. 
This fact, together with
Proposition~\ref{distortion} imply the next
proposition \cite{afryou98}.

\begin{prop}\label{density}
Let $\Om$ be a measurable subset of $[0,\bar{\ga})$ and denote
$\Om_l = \Om \cap [\ga_l,\ga_{l-1}]$. If the limit
$$
\lim_{l \rightarrow +\infty} m(\theta_l(\Om_l)) 
$$
exists and equals $\De$, then
$$
\lim_{\ga \searrow 0}
\frac{m\left( \Om \cap [0 , \ga)\right) }{\ga} = \De.
$$
\end{prop}

We remark that, with $\Omega_l$ satisfying the assumption in the above proposition,   
$m(\Omega_l) / \ga_l \sqrt{\ga_l}$ converges as $l \to \infty$.

Let $\Lc_{l,\theta}$  denote the 
local (first hit) map from $I^s_{g_l(\theta)}$ to $I^u_{g_l(\theta)}$ induced by $f_{g_l(\theta)}$. 
The convenience of using $\tau^s_\ga$ and $\tau^u_\ga$ as coordinates on 
$I^s_\ga$ and $I^u_\ga$ is seen in the following proposition.

\begin{prop}\label{localmap}
For each $l \ge l_0$ and each $\theta \in [0,1]$
\begin{equation}\label{eq_localmap}
\Lc_{l,\theta}
 = (\tau^u_{g_l(\theta)})^{-1} \circ R_{-\theta} \circ \tau^s_{g_l(\theta)},
\end{equation} 
where $R_\theta$ denotes a rigid rotation by angle $\theta$.
\end{prop}

\noindent {\sc Proof.}
This follows from Proposition~\ref{embedding} and the definitions 
of $\tau_{g_l(\theta)}^s$ and $\tau_{g_l(\theta)}^u$
as the time variables for the embedding flow for $f^q_{g_l(\theta)}$.
\qed \\

In addition to $\tau^s_{\ga} : I^s_{\ga} \rightarrow [0,1]$ defined 
above, it is also convenient to define an extension,
$\taub^s_{\ga}$, by
$$
      \phi_{\ga}^{\taub^s_{\ga}(x)}(d) = x.
$$
Thus $\taub^s_0:[d,a) \rightarrow [0,\infty)$ and 
$\taub^s_{\ga}:[d,e] \rightarrow [0,\infty)$ for $\ga >0$.
Furthermore, we will denote
$$
   \tau^s_{\ga}(x) = \taub^s_{\ga}(x) \mod 1
$$
wherever $\taub^s_{\ga}$ is defined.
Similarly, $\taub_\ga^u$ is given by
$$
      \phi_{\ga}^{\taub^u_{\ga}(x)}(e) = x.
$$

Let $c_i(\ga) = f^i_\ga(c)$, $i \ge 1$. 
It follows from the assumptions that there is an integer $j$ such that
$f^j_0(c) \in I^s_0$. By choosing $d$ so that
$f^j_0(c)$ is in the interior of $I^s_0$ we have
that $f^j_\ga(c) \in I^s_\ga$ for all $\ga$ small.
Denote this point by $c^s(\ga)$.
It then follows that for any $l$, either
$$
   c_{j + l q}(\ga) \in I^u_{\ga}
\qquad \text{ or } \qquad
   c_{j + (l+1) q}(\ga) \in I^u_{\ga}.
$$
Denote this intersection of $\{c_i(\ga)\}$ with $I^u_\ga$
by $c^u(\ga)$.
Note that for a fixed $l$ the function $\theta \mapsto c^u(g_l(\theta))$
will have a jump discontinuity at which the value will
jump from one endpoint of $I^u_{\ga}$ to the other. Denote
by $\theta^\sharp_l$ the point at which the discontinuity
takes place.

\begin{lem}
There exists a limit 
$$
    \lim_{l \rightarrow \infty} 
           \theta^\sharp_l \rightarrow \theta^\sharp_\infty.
$$
As $l \rightarrow \infty$ the sequence of maps 
$\theta \mapsto c^u(g_l(\theta))$ converges in the 
$C^r$ topology on compact sets not containing 
$\theta^\sharp_\infty$.
\end{lem}

\noindent {\sc Proof.}
This easily follows from Proposition~\ref{embedding} and 
the definition of $g_l(\theta)$.
\qed \\

Specifically, we will make use of the implication that
the derivatives of $c^u(g_l(\theta))$ with respect to $\theta$
converge uniformly as $l \rightarrow \infty$ for
$\theta$ in compact intervals away from the discontinuity.
Later, this will allow us to make estimates of derivatives
along $\{c_i(\ga)\}$ which are
uniform in $l$. 


\section{The saddle-node: global analysis}
\label{sec_mather}

In this section, 
we apply the local analysis near the saddle-node from the previous section
to obtain expressions for global return maps. 
We study the occurrence of Misiurewicz maps for small parameter values.

\subsection{The Mather invariant and return maps}

We begin with a simple but, useful lemma.

\begin{lem}\label{lem_over}
$W^u(a) = [f^2_0(c) , f_0(c)]$.
\end{lem}

\noindent {\sc Proof.}
Let $N$ be the interval bounded by $a$ such that
$f^q_0 (N) = N$.
Consider 
$$
M =  \bigcup_{i\ge 0} f_0^{iq} (W^u_{loc} (a)).
$$
Note that $M$ is an invariant set.
We have $f^q_0 (N \cup M) = N \cup M$.
If $M$ does not cover $[f^2_0(c) , f_0(c)]$,
then for some $0\le i < q$, $f^i_0   (N\cup M)$ is disjoint from $N \cup M$.
This contradicts the assumption
that $f_0$ is once renormalizable. 
\qed\\

We will use the freedom in the choice of $e$ and take  
$e$ so that $c$ is not in $\Oc(e)$.
By this choice we have that
$c$ will be in the interior of $f^i_0(I^u_0)$ for some $i$.

For $\ga = 0$, let $\bar{\Gc}$
be the first hit map from $I^u_0$ to $W_{loc}^s(a) = [d,a)$.
Define $\bar{\Mb}: [0,1] \rightarrow \bR$ by 
$$
    \bar{\Mb} = \taub^s_0 \circ \bar{\Gc} \circ (\tau^u_0)^{-1},
$$ 
and  define $\Mb: [0,1] \rightarrow [0,1]$ by $\Mb = \bar{\Mb}  \mod 1$. 
By identifying the endpoints of $[0,1]$ we may consider $\bar{\Mb}$
as a map from a subset $D(\bar{\Mb})$ of $\bS^1$ into $\bR$ and $\Mb$ as a map from 
$D(\Mb) = D(\bar{\Mb})$ into $\bS^1$.
Following \cite{yoc95}, we call $\bar{\Mb}$ the {\em Mather invariant} for $f_0$.
One may show that $\Mb$ is a modulus of 
smooth conjugation, in other words, it is invariant under differentiable
changes of variables.

For each $\tau$ in the domain, $\Dc(\Mb)$, of $\Mb$ and $x \in I^u_0$ 
such that $\tau = \tau^u_0(x)$ associate the integer
$$
\Rb(\tau) = \min_{x_i \in \Oc(x) \cap (a,f^q_0(e))} \lfloor \taub^u_0(x_i) \rfloor.
$$ 
Some $x \in I^u_0$ will return (perhaps many times) to $(a,f^q_0(e)]$ before eventually being
trapped in $W^s(a) = [d,a)$ and
the function $\Rb :[0,1] \rightarrow \bZ^-$  measures the depth of the
deepest return. 

Finally, to each $\tau$ in the domain of $\Mb$ and the corresponding point $x \in I^u_0$
associate a positive integer
$$
\Nb(\tau) = \min\{i | f^i_0(x) \in W^s(a)\}.
$$

Given $j$, denote by $V(j)$ the subset of
$\Dc(\Mb)$ defined by
$$
V(j) = 
  \{\tau \in \Dc(\Mb) :  \bar{\Mb}(\tau) < j, 
         |\Rb(\tau)| < j     
         \text{ and }  \Nb(\tau) < j \}.
$$ 
The points in $(\tau^u_0)^{-1}(V(j)) \in I^u_0$ are those whose forward
orbits enter $W^s_{loc}(a)$ in a bounded number iterations 
and which do not come too close to $a$ in the process, 
either by re-entering $W^u_{loc}(a)$ too close to a
or by landing in $W^s_{loc}(a)$ too close to $a$.
Since the forward orbit of almost every $x \in W^u_{loc}(a)$ has $\Oc(a)$ as its 
omega limit set, the measure $m(V(j))$ may be made close to $1$ 
by taking $j$ to be large.

We will let $\bar{E} = \cup_{i=0}^{q-1} f^i_0(E)$, so that
$\bar{E}$ is a neighborhood of the orbit of $a$ for $\ga = 0$.
Denote by $e_{-i} = \left( \Res{   f^{-iq}_{\ga} }{E} \right)^{-1}(e)$ and
\begin{equation}\label{Ei}
E^{-i}_{\ga} = [e_{-i-1},e_{-i}).
\end{equation}
Note that orbits associated with $V(i-1)$ do not 
reenter $E^{-i}_0$ before they hit $[d,a)$.
Also note that for $i$ fixed, the intervals $E^{-i}_\ga$ approach
$E^{-i}_0$ as $\ga \rightarrow 0$.

Consider the first return map, 
$\ka_{i,\ga}$, of the interval $E^{-i}_{\ga}$ and
let $\tilde{\ka}_{i,\ga}$ be the normalized map given by
$$
 \tilde{\ka}_{i,\ga} = 
      \tau^u_{\ga}|_{E^{-i}} \circ \ka_{i,\ga} 
            \circ (\tau^u_{\ga}|_{E^{-i}})^{-1}.
$$
Identifying the endpoints of $[0,1]$, we may consider
$\tilde{\ka}_{i,\ga}$ as a map on the circle.

\begin{prop}\label{approx}
Given any $i$, 
$$
\lim_{l \rightarrow \infty}  \Bigl| \tilde{\ka}_{i,g_l(\theta)}|_{V(i-1)} 
        -  R_{-\theta} \circ \Mb|_{V(i-1)} \Bigr|_{C^r} \rightarrow 0
$$
for each $0\le \theta < 1$.
\end{prop} 

\noindent {\sc Proof.}
Denote $d_i = f^{iq}_{\ga}(d)$ and $D^i = [d_i,d_{i+1})$.
Let $\Gc_{i,\ga}$ denote the global (first hit) map from 
$E^{-i}$ to $D^{i}$induced by $f_{\ga}$. 
Let $\Lc_{i,\ga}$ denote the local map.
Note that $\ka_{i,\ga}$ is not equal to 
$\Lc_{i,\ga} \circ \Gc_{i,\ga}$ since some points in $E^{-i}$ 
will return to $E^{-i}$ before hitting $D^{i}$ (by landing 
in $[d_{i+1},e_{-i})$). However, the two maps do agree when 
restricted to $(\Res{\tau^u}{E^{-i}})^{-1}(V(i-1))$ since 
those points will in fact hit $D^{i}$ before returning to $E^{-i}$.
Since we are only considering a finite number of iterations 
from $E^{-i}$ to $D^{i}$, it follows from the construction 
of $\Mb$ that $\tau^s_{i,\ga} \circ \Gc_{i,\ga} \circ (\tau^u_{i,\ga})^{-1}$ 
converges to $\Mb$ on the restricted set $V(i-1)$.
Proposition~\ref{localmap} then implies that, 
\begin{equation}
\tau^u_{\ga} \circ \Lc_{j,\ga} \circ \Gc_{i,\ga} 
          \circ (\tau^u_{\ga}|_{E^{-i}})^{-1}
\end{equation} 
converges to $R_{-\theta} \circ \Mb$ on $V(i-1)$.
\qed 


\subsection{Misiurewicz maps}
\label{sec_misiur}

In this section we will identify parameter values $\theta$
for which $f_{g_l(\theta)}$ maps the critical point $c$ onto some
repelling hyperbolic periodic point and thus $c$ does not
return to a neighborhood of itself, i.e. $f_{g_l(\theta)}$ satisfies the
{\em Misiurewicz condition}. There are two ways that
hyperbolic periodic points can occur: as periodic
points whose orbits pass through $\bar{E}$, or as
continuations of periodic points for $\ga =0$ 
(outside of $\bar{E})$. 
In the former case, the periodic orbits exist for parameter
values within subintervals of $[\ga_{l+1},\ga_l]$ for each $l$.
In the latter case
the periodic orbits exist for all $\ga >0$ sufficiently small.
We confine ourselves to a discussion of this case. 
We show that
parameters for which $c$ is mapped onto
the periodic orbits under consideration occur in decreasing
sequences $\{\ga_l^*\}$ of values approaching $0$ as 
$l \rightarrow \infty$. Later we will show that 
such a sequence typically gives rise to the parameter sets predicted 
in Theorem~\ref{existence}.

Note that the assumptions on $f_0$ imply
that the nonwandering set, $L_0$, of $f_0$ restricted to 
$[0,1] \setminus \bar{E}$ is a hyperbolic set 
(see Theorem III.5.1  in \cite{melstr93}). 
If a point $x \in I^u_0$ is mapped onto $L_0$ then 
$\Mb$ will have a discontinuity at $\tau^u_0(x)$.
By Lemma~\ref{lem_over}, for any $y^* \in L_0$ there is a point $x^* \in I^u_0$ which
is mapped onto $y^*$ by $f_0$.
Recall that $c_j(\ga) \in I^s_\ga$ is the first point in the orbit of $c$ that hits
$I^s_\ga$. Let $\tau^c_0 = \tau^s_0 (c_j (0))$.

\begin{lem}\label{Mis2}
Suppose that $y^*_0$ is a periodic point contained in the nonwandering set
$L_0$ of $[0,1] \setminus \bar{E}$. Denote by $y^*_\ga$ the continuation of 
$y^*_0$ for $\ga \ge 0$. 
Let $\tau^*$ be a point of discontinuity of $\Mb$ corresponding to $y^*_0$ 
and suppose that $\theta^*$ is such that $\tau^c_0$ 
is mapped onto $\tau^*$ 
by $(R_{-\theta^*} \circ \Mb)^k \circ R_{-\theta^*}$ for some $k \ge 0$. 
Suppose further that 
\begin{equation}\label{partialMb}
 \frac{\partial}{\partial \theta}\Big|_{\theta = \theta^*} 
    \left( (R_{-\theta} \circ \Mb)^k \circ R_{-\theta}( \tau^c_0 ) 
                           - \tau^* \right) \ne 0.
\end{equation}
Then 
there is a sequence of parameter values 
$\{\ga_l^*\}$, $\ga_l^* \in [\ga_{l+1},\ga_l]$, such that 
$c$ is mapped onto $y^*_{\ga_l^*}$ by iterates of $f_{\ga_l^*}$.
\end{lem}

\noindent {\sc Proof.}
The point  
$c_j(\ga) = f_{\ga}^j(c)$ is in  $I^s_{\ga}$
for all small $\ga$, and,  
$\tau^c_{\ga} = \tau^s_{\ga}(c_j(\ga))$ approaches 
$\tau^c_0$ as $\ga \searrow 0$. Thus, 
$c^u(\ga) = f^{j+ql}_{\ga}(c) \in I^u_{\ga}$. 
Let $x^*$ be the point in $I^u_0$ which is mapped
onto $y^*_0$ by an iterate of $f_0$.
If $k = 0$, then $c^u(g_l(\theta^*))$ approaches $x^*$ as
$l \rightarrow \infty$ and (\ref{partialMb}) implies the result.
If $k>0$ then it follows there is an $m$ such that 
$(R_{-\theta^*} \circ \Mb)^i(\tau^*) \in V(m-1)$
for all $0\le i \le k-1$.
Applying Proposition~\ref{approx} to $D^m$ and $E^{-m}$ then implies 
that $c^u(l,\theta)$ is
mapped by the $k$-return map into $E^{-m}_{g_l(\theta^*)}$ arbitrarily
close the $mq$-th preimage of $x^*$ for $l$ large. The result then follows
from  (\ref{partialMb}).
\qed\\

The following corollary provides the parameter values which will be used in
the proof of Theorem~\ref{existence}.

\begin{prop}\label{oml}
Given any $y^* \in L_0$, there exists integers $m >0$ and $l_0 >0$
and a sequence of parameter values $\ga_l^* \in [\ga_{l+1},\ga_l]$,
$l\ge l_0$ such that
$$
    f^{lq+m}_{\ga_l^*}(c) = y^*_{\ga_l^*}.
$$
Further, there is a number $\theta^* \in (0,1)$ such that 
$\theta_l(\ga_l^*) \rightarrow \theta^*$ as $l \rightarrow \infty$.
\end{prop}

\noindent {\sc Proof.}
For each $y^* \in L_0$ there is a point $x^* \in I^u_0$ such
that $x^*$ is mapped to $y^*$ by an iterate of $f_0$. Let $\tau^* = \tau^u_0(x^*)$.
If we let $\theta^* = 1 - \tau^* + \tau^c$ then 
$R_{-\theta^*}(\tau^c) = \tau^*$ and 
$$
 \frac{\partial}{\partial \theta}\Big|_{\theta = \theta^*} 
    \left( R_{-\theta}( \tau^c_0 ) 
                           - \tau^* \right) = 1.
$$
Then Lemma~\ref{Mis2} implies
the result.
\qed




\section{Parameter values with bounded recurrence}\label{sec_br}

We start the proof of Theorem~\ref{existence}.
In this section we construct a set $\Omega$ 
with positive density at $\ga = 0$, so that
$f_\ga$ for $\ga \in \Omega$ has bounded recurrence (see Definition~\ref{def_br};
the definition is in terms of an induced map which is defined below).
We deduce that $f_\ga$ for $\ga \in \Omega$ has
an absolutely continuous invariant measure.

Considering $f_\ga$ for $\ga \in [\ga_{l+1},\ga_l]$, we
reparameterize and use the parameter $\theta \in [0,1]$
defined by $g_l(\theta) = \ga$, see Section~\ref{sec_snlocal}.
Denote the resulting family of unimodal maps by $f_{l,\theta}$.
By Proposition~\ref{oml}, there exists a converging 
sequence of parameter values
$\theta_l^* = \theta_l(\ga_l^*)$ 
so that $f^{q l + m}_{l,\theta_l^*} (c)$ is a hyperbolic
periodic point.
The map $f_{l,\theta_l^*}$ is a Misiurewicz map and by an extension
of Jakobson's Theorem, see Theorem V.6.1 in \cite{melstr93}, $\theta_l^*$ is a 
Lebesgue density point of a set of parameters $\Theta_l$ for which
$f_{l,\theta}$ supports an absolutely continuous invariant measure.
Theorem~\ref{existence} can be proved by establishing that 
the measure of $\Theta_l$ is bounded away from 0 uniformly in $l$.
Proposition~\ref{density} guarantees that the union 
$\cup_l g_l(\Theta_l)$ has positive Lebesgue measure and 
positive density at $\ga = 0$.
Combined with the statement that $f_{l,\theta}$, $\theta \in \Theta_l$, 
admits an absolutely continuous invariant measure, this proves Theorem~\ref{existence}. 

In the following section we define induced maps. 
The actual construction of the sets
$\Theta_l$ is described in Section~\ref{sec_indconstructions}.

\subsection{Induced maps}\label{sec_ind}

As before, $a$ denotes a saddle-node periodic point of $f_0$,
of period $q$.
We may assume that $a$ is nearest to the critical point
$c$ of all points in ${\cal O}(a)$,
so that $f^q_0$ is a homeomorphism on $(c,a)$.

We may suppose that there is a periodic point $z_0 \in L_0$
and an integer $j$ such that $f_0^j(e) = z_0$. Let $z_\ga$
be the periodic continuation of $z_0$ and $e_\ga$ the 
continuation of $e_0$ such that $f_0^j(e_\ga) = z_\ga$.
We will suppress the $\ga$ dependence of $e$ in the notation.

Let $E^{-i}_{l,\theta}$ be the fundamental intervals given
by $(e_{-i},e_{-i+1})$. Denote 
\begin{eqnarray*}
\Et_{l,\theta} & = & \cup_{i=0}^l E^{-i}_{l,\theta} \subset E.
\end{eqnarray*}
Note that $[d_1,e) \subset \Et_{l,\theta}$.
Before proceeding with the definition and construction of $\Theta_l$, 
we introduce an induced map $\ft_{l,\theta}$.

Namely, 
$\ft_{l,\theta}$ will be defined from $f_{l,\theta}$ in the following 
way, $\ft_{l,\theta}$ will equal $f_{l,\theta}$, except for points 
in $\Et_{l,\theta}$. Points which fall in this set
will be mapped ahead to the image of $I^u_{l,\theta} = E^0_{l,\theta}$.
That is,
\begin{equation}\label{tildef1}
\ft_{l,\theta}(x) = \left\{ \begin{array}{ccc}
          f_{l,\theta}^{iq + 1}(x), 
         & \quad & \text{if } x \in E^{-i}_{l,\theta} \\
        f_{l,\theta}(x), & \quad & \text{otherwise}.
                     \end{array} \right.
\end{equation}
Thus points in $\Et$ are mapped by the first hit map to the interval
$f_{l,\theta}(I^u_{l,\theta})$. For $\ga \in (\ga_{l+1},\ga_l)$, 
$\ft_{l,\theta}$ will have $l+1$ discontinuities.
Recall that for $\ga >0$ and $x \in [d,e)$ the function
$\tau^s_{l,\theta}$ is defined by $\tau^s(x) = \taub^s(x) \mod 1$ where
$$
   \phi^{\taub^s_{l,\theta}(x)}_{l,\theta}(d) = x.
$$
Observe that $\ft_{l,\theta}$ equals
\begin{equation}\label{tildef}
\ft_{l,\theta}(x) = \left\{ \begin{array}{ccc}
     f_{l,\theta} \circ (\tau^u_{g_l(\theta)})^{-1} \circ \tau^s_{g_l(\theta)}(x), 
         & \quad & \text{if } x \in \Et_{l,\theta} \\
        f_{l,\theta}(x), & \quad & \text{otherwise}.
                     \end{array} \right.
\end{equation}

Write
$$
\Et_{\infty,\theta} = 
\left( 
(\taub_0^s)^{-1} \circ R_\theta \circ \taub_0^u (\Res{f_{l,\theta}}{E})^{-q} (e) , e
\right).
$$
Let $\ft_{\infty,\theta}$ denote the map defined from $f_0$ in the following way.
If $x \notin \Et_{\infty,\theta}$ then $\ft_{\infty,\theta}(x) = f_0(x)$. If $x \in (a,e)$
then $\ft_{\infty,\theta}(x)$ is given be the first hit map from $(a,e)$
to $f_0(I^u_0)$ induced by $f_0$. If $x \in (d,a)\cap \Et_{\infty,\theta}$ then let
$\tau_0^s(x)$ be defined as above and let 
\begin{equation}\label{tildefinfty}
\ft_{\infty,\theta}(x) = f_0 \circ (\tau_0^u)^{-1} \circ
        R_{-\theta} \circ \tau_0^s(x) .
\end{equation}
In other words, if $x \in (d,a) \cap  \Et_{\infty,\theta}$ then $x$ is mapped to
$f_0(I^u_0)$ using the time coordinate shifted by $-\theta$.
Note that $\ft_{\infty,\theta}(x)$ will have countably many discontinuities
accumulating at $a$ from both sides. In $(a,e)$ these occur at
$\{e_{-i}\}_{i=0}^\infty$. In $(d,a)$ these occur at those
$x$ for which 
$$
  \tau_0^s(x) = \theta.
$$

\begin{prop}\label{approx2}
For each $\theta$ the sequence of maps 
$\{\ft_{l,\theta}\}_{l=l_0}^\infty$ converges
to the map $\ft_{\infty,\theta}$, uniformly on compact sets 
which are disjoint from the discontinuities of $\ft_{\infty,\theta}$.
\end{prop}

\noindent {\sc Proof.}
Given that $\ft_{l,\theta}$ satisfies (\ref{tildef}) and
$\ft_{\infty,\theta}$ is defined by (\ref{tildefinfty}), this 
follows from Proposition~\ref{localmap}.
\qed\\

Observe that, with $\theta^*$ and $m$ as in Proposition~\ref{oml},
$\ft_{\infty,\theta^*}^{m}(c) = y^*$.

\subsection{Iterating intervals}

In the following we will consider iterates of intervals.
The map $\tilde{f}_{l,\theta}$, $l_0 \le l < \infty$, is 
discontinuous along backward iterates of $e$ 
in $E$, so that the intervals in the image of an interval
under $\tilde{f}_{l,\theta}$ might be arbitrarily small, regardless
of the size of the original interval.
We therefore slightly adjust the definition of $\tilde{f}_{l,\theta}$
to avoid this problem. Consider an interval $I \subset [0,1]$.
For $0\le i\le l$, denote $I^{-i} = I \cap E^{-i}_{l,\theta}$
and let $I^{-l-1}$ be the component of
$I \backslash E$ adjacent to $E^{-l}_{l,\theta}$ and
$I^1$ the component 
of $I \backslash E$ adjacent to $E^0_{l,\theta}$.
This yields a partition $\{ I^i\}_{-l-1 \le i \le 1}$ of $I$
(elements of the partition can be empty).
If the leftmost or rightmost nonempty intervals of this partition
do not contain a fundamental domain 
$E^{-i}_{l,\theta}$,
$-l-1 \le -i \le 1$, join them to the adjacent intervals.
Note that if $I$ is partitioned into two elements $\{I^{-i},I^{-i+1}\}$ neither of which
is a fundamental domain, this leaves a choice in coding the resulting
interval after $I^{-i}$ or $I^{-i+1}$.
This way an interval $I$ that covers
one or more fundamental domains
is partitioned 
into subintervals which are at least as large
as a fundamental domain.
Given $x \in I^i \subset I$, 
define
\begin{eqnarray*}
\breve{f}_{l,\theta} (x;I)  & = &
\left\{ 
\begin{array}{ll}
f_{l,\theta} (x), & \mbox{if } i = -1 \mbox{ or }i = l + 1,
\\
f_{l,\theta}^{iq + 1} (x),  & \mbox{if }  0\le i \le l.
\end{array}
\right.
\end{eqnarray*}
Note that as long as $I$ does not cover a fundamental domain
in $E$, $\breve{f}_{l,\theta} (\cdot;I)$ equals some fixed 
iterate of $f_{l,\theta}$. Also note that 
$\breve{f}_{l,\theta} (I;I)$ consists of at most two
components. We further remark that, if 
$\breve{f}_{l,\theta}^j$ maps $I$ homeomorphically onto 
$E^0_{l,\theta}$, then there is a fixed number $N$ of 
iterates after which $\breve{f}_{l,\theta}^{j+N} (I;I)$
contains $c$ in its interior.

Define maps $F_l$ and $\tilde{F}_l$ by
\begin{eqnarray*}
F_l(x,\theta) &  = & (f_{l,\theta}(x),\theta),
\\
\tilde{F}_l(x,\theta) & = & (\tilde{f}_{l,\theta}(x),\theta).
\end{eqnarray*}
We would like to consider images $\tilde{F}_l^i (0,\theta)$.
As for single maps
we come across the difficulty 
that $\ft_{l,\theta}$ is discontinuous along the backward orbit 
of $e$ in $E$.
We consider the set $T = (\theta^*-\ep,\theta^*+\ep)$.
By a fundamental strip we mean
a set $\{ ( E^{-i}_{l,\theta}, \theta)\}$, $\theta \in T$.
Consider a curve $C = \{x(\theta),\theta)\}$ that 
projects injectively to $[0,1]$ by the projection
$(x,\theta) \to x$.
Define $C^{-i}$ to be the intersection of  $C$  with the fundamental strip
$\{ ( E^{-i}_{l,\theta}, \theta)\}$, $\theta \in T$,
and let
$C^{-l-1}$ be the connected component of $C$ that
is adjacent to $C^{-l}$ and
$C^{1}$ the connected component
adjacent to $C^0$.
This defines a partition $\{ C^{-l-1}, \ldots, C^1 \}$
of $C$ with possibly empty elements.
If the leftmost or rightmost nonempty element of this partition
does not cross a fundamental strip,
join it to the adjacent element.
This way a partition of a curve $C$ that 
crosses at least one fundamental strip 
is obtained all of whose elements cross a fundamental strip.
Define 
\begin{eqnarray}\label{breveF}
\breve{F}_l (x,\theta) & = &
\left\{
\begin{array}{ll}
F_l(x,\theta), & \mbox{if } (x,\theta) \in C^{-l-1} \cup C^1,
\\
F^{iq+1}_l (x,\theta), & \mbox{if }
(x,\theta) \in C^{-i}, 0\leq i \le l.
\end{array}
\right.
\end{eqnarray}

\subsection{Inductive constructions}\label{sec_indconstructions}

In this section the actual construction of the set $\Theta_l$ 
is described. We make use of the induced maps introduced in 
Section~\ref{sec_ind} and for the rest we closely follow \cite{luz00}.

Consider $\theta$ near $\theta_l^*$ for a fixed value of $l$.
Write 
\begin{eqnarray*}
\ct_i({l,\theta}) & = & \ft_{l,\theta}^i (c).
\end{eqnarray*}

\begin{defn}\label{def_br}
For $\delta >0, \alpha > 0$, we say that $\ft_{l,\theta}$ 
satisfies the bounded recurrence condition $(\BR)_n = 
(\BR)_n (\alpha ,\delta)$ if for all positive integers $k \le n$,
\begin{equation}\label{br_n}
\prod_{\substack{ \ct_i({l,\theta}) \in (c-\delta,c+\delta),\\  
0\le i \le k}} |\ct_i({l,\theta}) - c|   \ge  e^{-\alpha k}.
\end{equation}
We say that $\ft_{l,\theta}$ satisfies $(\BR)$ if it satisfies $(\BR)_n$
for all $n$.
\end{defn}

Similarly, we say that $f_{l,\theta}$ satisfies 
$(\BR)_n$ if it satisfies the equivalent 
condition on $\{c_i(l,\theta)\}$. It is easy to see that if 
$\ft_{l,\theta}$ satisfies $(\BR)_n$ then so does $f_{\l,\theta}$,
and vice versa.
The next proposition is the main result of Section~\ref{sec_br} 
and will be shown to imply Theorem~\ref{existence}.

\begin{prop}\label{jakobson}
For every $\alpha > 0$ there exists $\delta >0$, so that
there is a set $\Theta_l$ of parameter values
for which $\ft_{l,\theta}$, $\theta \in \Theta_l$, satisfies $(\BR)$ and
\begin{eqnarray*}
\liminf_{\varepsilon\searrow 0} 
\frac{|\Theta_l \cap (\theta_l^* - \varepsilon,\theta_l^* + \varepsilon  )|}
     { 2 \varepsilon } & > & \sigma,
\end{eqnarray*}
for some $\sigma >0$, uniformly in $l$.
\end{prop}

By Proposition~\ref{density}, 
$\cup_l g_l(\Theta_l)$ has positive density at $\ga = 0$.
Theorem~\ref{existence} is shown by establishing that 
$f_{l,\theta}$ possesses an absolutely continuous invariant measure
for $\theta \in \Theta_l$.
This is done more formally in Section~\ref{sec_proof},
but the idea is straightforward: one shows that 
$\ft_{l,\theta}$, $\theta \in \Theta_l$,
is a Collet-Eckmann map (meaning that
$|D\ft^i_\ga (f_\ga(c))| \ge K \lambda^i$ for some $K>0$, $\lambda >1$).
This implies that also $f_{l,\theta}$, $\theta \in \Theta_l$,
is a Collet-Eckmann map, albeit with weaker expansion.
Collet-Eckmann maps are known to admit absolutely continuous
invariant measures.

To prove Proposition~\ref{jakobson}, we adapt the reasoning
in \cite{luz00}, where Jakobson's result is proven
using a variant of Benedicks-Carleson's proof. 
Below we will describe the constructions and show key estimates
and computations that differ from those in \cite{luz00}.
After some preliminary work, the reasoning will follow
\cite{luz00} closely. 
We will therefore be able to simply refer
to \cite{luz00} for several of the more technical parts
of the proof.

Note that we must consider a sequence of families $\ft_{l,\theta}$ 
that tends to the family $\ft_{\infty,\theta}$ as $l \to \infty$, 
see Proposition~\ref{approx2}. It is instructive to compare with 
\cite{diarocvia96} where the analysis leads to the study of a 
sequence of families $g_{l,\theta}$ of smooth unimodal maps 
(renormalizations of $f$), 
converging to the logistic family as $l \to \infty$. As they 
remarked, there are uniform lower bounds on the measure of the 
set of parameter values for which absolutely continuous invariant 
measures occur, for families from a neighborhood of the logistic 
family (and so for all large enough values of $l$).
A similar situation was considered in \cite{pumrod97}.
Following the reasoning below one can construct a set of parameter 
values for which $\ft_{\infty,\theta}$ has bounded recurrence 
and show that this set has positive measure.
In fact, the reasoning below constructs such sets, and
uniformly bounds their measures, for families nearby
$\ft_{\infty,\theta}$, in particular for the families $\ft_{l,\theta}$
with $l$ large. We concentrate on the families $\ft_{l,\theta}$
and formulate the results in terms of this sequence of families.

We now start the constructions.
For a positive integer $r$, let $I_r = [c + e^{-r},c + e^{-r+1})$ and
$I_{-r} = (c - e^{-r+1},c - e^{-r}]$.
Let $\iota$ be a small positive number.
Given $\delta >0$, write $r_\delta = -\ln \delta$ and
$r_{\delta^+} = - \iota \ln \delta$.
We can suppose that $r_\delta$ and $r_{\delta^+}$ are integers.
Let
\begin{eqnarray*}
\Delta & = & \{c\} \cup \bigcup_{|r| \ge r_\delta +1} I_r,
\\
\Delta^+ & = & \{c\} \cup \bigcup_{|r| \ge r_{\delta^+} +1} I_r.
\end{eqnarray*}
Subdividing each interval $I_r$ into $r^2$ subintervals, $I_{r,m}$, of equal length
provides partitions ${\cal I}$ of $\Delta$ and ${\cal I}^+$ of $\Delta^+$.

Given $x \in \Delta^{+}$, write $\eta_0 = (c,x)$ (or $(x,c)$) and let
$\eta_j =  \breve{f}^j_{l,\theta} (\eta_0;\eta_0)$.
Define the {\em binding period} of $x$ as 
\begin{eqnarray}\label{binding1}
q_l(x,\theta) & = & \sup \{ m \in \bN: 
 \left| \eta_j \right| 
\le e^{-2\alpha j} \mbox{ for all } 0 \le j \le m-1 \}
\end{eqnarray}.
Suppose $\tilde{c}_k(l,\theta) \in \Delta^+$ and define the binding period 
associated with $\tilde{c}_k(l,\theta)$ as
$$
   p(l,\theta,k) = q_l(\tilde{c}_k(l,\theta),\theta).
$$
For a fixed $l$, let $\omega$ be a subinterval of 
$(\theta_l^* - \epsilon , \theta_l^* + \epsilon)$.
Denote 
\begin{eqnarray*}
\omega_i & = & \Pi \circ \breve{F}^i_l (c,\omega),
\end{eqnarray*} 
where $\Pi$ is the projection $\Pi(x,\theta) = x$.
If $\omega_k$ intersects $\Delta$, $0\le k \le n-1$, then $k$ is called a {\em return time} for $\omega$.
Define the binding period associated with a return time $k$ of 
a parameter interval $\omega$ as
\begin{equation}\label{binding2}
p(\omega,k)  =  
      \min_{\substack{ \theta \in \omega, \\
             \tilde{c}_k(l,\theta) \in \Delta^+ }}
              p(l,\theta,k).
\end{equation}

For each $l$ let $\Pc_l^{(0)}$ be the trivial partition 
$\{ (\theta_l^* - \epsilon , \theta_l^* + \epsilon)\}$ of
the parameter interval 
$\Theta_l^{(0)} = (\theta_l^* - \epsilon , \theta_l^* + \epsilon)$. 
Inductively we will define parameter sets $\Theta_l^{(n)}$
and partitions $\Pc_l^{(n)}$ thereof.
In order to define $\Pc_l^{(n)}$ given $\Pc_l^{(n-1)}$, we first 
construct a refinement $\hat{\Pc}_l^{(n)}$ of $\Pc_l^{(n-1)}$.

We say that a return of $\omega$ at time $k$ is a {\em bound return} if 
there is a return time $j <k$ of $\omega$ and $k \le j + p(\omega,j)$.
Let $\omega \in \Pc_l^{(n-1)}$.

\begin{description}

\item[{\bf Chopping times.}]
We say that $n$ is a chopping time 
for $\omega$ if
\begin{enumerate}
\item
$\omega_n$ intersects $\Delta^+$ in at least
three elements of ${\cal I}^+$, and
\item
$\omega_n$ is not a bound return for $\omega$.
\end{enumerate}

\item[{\bf Non-chopping times.}]
We say that $n$ is a non-chopping time 
for $\omega$ in all other cases, that is 
if one or more of the following occurs:
\begin{enumerate}
\item
$\omega_n \cap \Delta^+ = \emptyset$.
\item
$\omega_n$ is a bound return of $\omega$.
\item
$\omega_n$ intersects 
no more than two elements of the partition ${\cal I}^+$ of $\Delta^+$.
\end{enumerate}

\end{description}

In case $n$ is a non-chopping time for $\omega$, we let
$\omega \in \hat{\Pc_l}^{(n)}$.
If $n$ is a chopping time for $\omega$ we partition $\omega$
as follows.
Write
$\omega_n \cap \Delta^+ = \cup_{m} \omega^{m}_n$, so that each
$\omega_n^{m}$ fully contains one and at most one
element of ${\cal I}^+$.
If $\omega_n \backslash \Delta^+$ contains an interval of length less then
$\delta^\iota$, we include this interval in the adjacent interval
of ${\cal I}^+$.
Otherwise an interval of $\omega_n \backslash \Delta^+$ is an
element of the partition of $\omega_n$.
Write the resulting partition of $\omega_n$ as 
$\omega_n = \cup_{m} \omega^{m}_n$. 
There is a corresponding partition 
$\cup_{m} \omega^{m}$ of $\omega$,
given by
$ \Pi\circ \breve{F}^n_l (c,\omega^m) = \omega^{m}_n$.
Let each element of this partition be an element of $\hat{\Pc_l}^{(n)}$.
Note that an element of $\hat{\Pc_l}^{(n)}$ partitioning $\omega$ need
not be connected, but can consist of several intervals.
This may happen if $\omega_j$ intersects $\Et_{l,\theta}$ in at least two
fundamental domains, for some $j<n$.

Let $\omega \in \Pc_l^{(n-1)}$ and consider 
$\nu \in \Res{ \hat{\Pc_l}^{(n)} }{\omega}$.
We speak of a bound, essential or inessential return time or an escape time
for $\nu$ in the following situations. \\

\noindent
{\em bound return time:}
The interval $\nu_n$ intersects $\Delta$ and
$n$ is a bound return time for $\omega$.\\

\noindent
{\em inessential return time:}
The interval $\nu_n$ intersects $\Delta$ and
$n$ is a non-chopping time for $\omega$ that is
not a bound return time, but
$\omega_n$ intersects at most two elements
of the partition ${\cal I}^+$.\\

\noindent
{\em essential return time:}
The interval $\nu_n$ intersects $\Delta$ and
$n$ is a chopping time for $\omega$.\\

\noindent 
{\em escape time:}
The return time $n$ is a chopping time for $\omega$, but
$\nu_n$ does not intersect $\Delta$. In this case we call
$\nu$ an escape component of $\omega$.\\

Any interval $\nu \in \hat{\Pc_l}^{(n)}$ belongs to a unique nested
sequence of intervals 
$$ \Theta_l^{(0)} \equiv \nu^{(0)} \supseteq \nu^{(1)} \supseteq \ldots
    \supseteq \nu^{(n-1)} \supseteq \nu^{(n)} \equiv \nu,
$$
where $\nu^{(k)} \in \Pc_l^{(k)}$ for $0 \le k < n$.
If $j$ is a chopping time for $\nu^{(j-1)}$, then $\nu^{(j)}$ is 
strictly contained in  $\nu^{(j-1)}$. Chopping times are either escape times
or essential return times.

The return depth of $\nu$ at time $k$
is defined if $\nu_k$ intersects $\Delta$, as
\begin{eqnarray*}
r & = & \max \{ |r|, \;\; \nu_k \cap I_r \neq \emptyset \}.
\end{eqnarray*}
Define functions
${\cal R}^{(n)} : \hat{\Pc}_l^{(n)} \to \bN$ and
${\cal E}^{(n)} : \hat{\Pc}_l^{(n)} \to \bN$ which associate to
$\nu \in \hat{\Pc}_l^{(n)}$ the sum of the return depths
and the sum of the essential return depths, over the first $n$
iterates $\Pi\circ \breve{F}^i_l(c, \nu)$.

Define
\begin{equation*} 
\Theta_l^{(n)}  =  \{ \nu \in \hat{\Pc_l}^{(n)}:
{\cal E}^{(n)} (\nu) \le \alpha n /10 \}
\end{equation*}
and
\begin{eqnarray*}
\Pc_l^{(n)} & = & \Res{\hat{\Pc_l}^{(n)}}{\Theta_l^{(n)}}.
\end{eqnarray*}
The sets
\begin{eqnarray*}
\Theta_l & = & \bigcap_n \Theta_l^{(n)}
\end{eqnarray*}
will be shown to satisfy the stated properties in 
Proposition~\ref{jakobson}.

\subsection{Expansion}\label{sec_exp}

Here expansion properties of the maps $f_{l,\theta}$, 
$\ft_{l,\theta}$ and $\breve{f}_{l,\theta}$ are discussed.
We relate expansion properties of these maps. In 
Section~\ref{sec_mane} we prove a Ma\~n\'e type result 
for  $\ft_{l,\theta}$, that is, we show that there is 
expansion along orbits outside a neighborhood of $c$.

\subsubsection{Expansion for induced maps}\label{sec_relate}

The relation between expansion along orbits of $f_{l,\theta}$, 
$\ft_{l,\theta}$ and $\breve{f}_{l,\theta}$ is discussed in the 
next two lemma's.

\begin{lem}\label{lemma_withouttilde}
If there exist $\tilde{C}>0$, $\tilde{\lambda} >1$ such that 
$| D\ft_{l,\theta}^n(x)  | \ge \tilde{C} \tilde{\lambda}^n$ for  
all $n >0$,
then there are $C_l > 0$ and $\lambda_l >1$, so that
$| Df^n_{l,\theta}(x) | \ge C_l \lambda^n_l$, for all $n>0$.
\end{lem}

\noindent {\sc Proof.}
Write $f_{l,\theta}^i(x) = f_{l,\theta}^{k(i)} \circ \tilde{f}_{l,\theta}^{m(i)}(x)$ 
with $k(i)$ minimal such nonnegative integer. Compute
$$
|Df^i_{l,\theta} (x)| \ge \min |Df^{k(i)}_{l,\theta}|  \tilde{C} 
( \tilde{\lambda}^{m(i)/i} )^i.
$$
Now $k(i)>0$ implies that the piece of orbit 
$\tilde{f}_{l,\theta}^{m(i)}(x),\ldots,
f_{l,\theta}^{k(i)}(\tilde{f}_{l,\theta}^{m(i)}(x))$
is in $\bar{E}$. Since $c \not\in \bar{E}$, the term $|Df_{l,\theta}|$ is 
bounded below in $\bar{E}$. Further, $k(i)$ is bounded above by 
$(l+1)q$ since any point in $\bar{E}$ is mapped outside of $\bar{E}$ in
$(l+1)q$ or fewer iterations. Therefore the quantity $|Df_{l,\theta}^{k(i)}|$ is bounded
from below by a constant $D_l$. 
Thus, 
$| Df^i_{l,\theta} (x) | \ge \tilde{C} D_l ( \tilde{\sigma}^{m(i)/i} )^i$.
We can let $C_l = \tilde{C} D_l$.
Since there is a minimum number of iterations
of $f_{l,\theta}$ needed for an orbit to enter $\bar{E}$ after leaving
$I^u$, and the number of consecutive iterations in $\bar{E}$ is
bounded above by $q(l+1)$, it follows that
the fraction $m(i)/i$ is bounded below by a constant $d_l >0$.
Hence, $\tilde{\sigma}^{m(i)/i}$ is strictly larger than some
number $\sigma_l > 1$.
\qed\\

Similarly one derives the following lemma relating 
expansion of $\breve{f}_{l,\theta}$ to expansion of $\ft_{l,\theta}$.

\begin{lem}\label{breve-tilde}
If $| D\breve{f}_{l,\theta}^n (x;I) | \ge \breve{C} \breve{\lambda}^n$ for some 
$\breve{C}>0$, $\breve{\lambda} > 1$,
then there are constants $\tilde{C} > 0$ and $\tilde{\lambda} >1$
so that
$| D\ft^n_{l,\theta} (x) | \ge \tilde{C} \tilde{\lambda}^n$. 
The converse statement holds as well.
\end{lem}

Given a subinterval $I\subset [0,1]$, 
write $\breve{f}^{j}_{l,\theta} (x;I) = \ft^{k(j)}_{l,\theta}(x)$.
Then $k(j) / j $ is bounded and bounded away from 0.
In fact, if $\Et_{l,\theta}$ is small, the quotient
$k(j)/j$ will be close to 1 for large values of  $j$, since 
there is then a large number of iterates
between any two passages through $\Et_{l,\theta}$.
Observe that for $\Et_{l,\theta}$ small, $\tilde{\lambda}$ in Lemma~\ref{breve-tilde} 
is close to $\breve{\lambda}$. 
Write $\ct_i (l,\theta) = \ft_{l,\theta}^i (c)$
and $\breve{c}_i (l,\theta) = \breve{f}_{l,\theta}^i (c;I)$.
Let $\breve{c}_j(l,\theta) = \ct_{k(j)} (l,\theta)$.
Then the bounded recurrence condition
\begin{equation*}
\prod_{\substack{ \ct_i({l,\theta}) \in (c-\delta,c+\delta), \\
     0\le i \le k }} |\ct_i({l,\theta})|  \ge  e^{-\alpha k},
\end{equation*}
where $k = k(j)$, translates into
\begin{equation*}
\prod_{\substack{ \breve{c}_i({l,\theta}) \in (c-\delta,c+\delta),\\
    0\le i \le j }} |\breve{c}_i({l,\theta})|  \ge  
e^{- C \alpha j},
\end{equation*}
for a constant $C$ which bounds $k(j) / j$ from below. 
If $\Et_{l,\theta}$ is small and $j$ is 
large, then $C$ is close to 1. Hence, when translating the bounded recurrence condition 
for $\ft_{l,\theta}$ into one for $\breve{f}_{l,\theta}$, an almost 
identical condition is obtained. A similar remark can be made 
for $\tilde{F}_l$ and $\breve{F}_l$.
   
\subsubsection{Ma\~n\'e's Theorem for induced maps}\label{sec_mane}

The following proposition is central and provides exponential
expansion of iterates of $\ft_{l,\theta}$ that stay away from 
the critical point $c$.
The estimates are uniform in $(l,\theta)$.
The proposition is modeled after
Theorem III.6.4 in \cite{melstr93}, which treats
families of smooth unimodal maps.

\begin{prop}\label{prop_exp}
There are constants $\varepsilon,\tilde{C} >0$ and 
$\tilde{\lambda} >1$, and a neighborhood $W$ of $c$ so that 
for any neighborhood $U$ of $c$ with $U \subset W$ and each 
large enough integer $l$, the following holds.
For each  $|\theta - \theta_l^*| < \varepsilon $,
if $\ft^j_{l,\theta} (x) \not \in U$ for $0\le j \le m-1$ and
$\ft^m_{l,\theta} (x) \in W$, then
\begin{eqnarray*}
\left| D\ft^m_{l,\theta} (x) \right| & \ge & \tilde{C} \tilde{\lambda}^m.
\end{eqnarray*}
If $\ft^j_{l,\theta} (x) \not \in U$ for $0\le j \le m-1$,
but not necessarily
$\ft^m_{l,\theta} \in W$,
then
\begin{eqnarray*}
\left| D\ft^m_{l,\theta} (x) \right| & \ge & 
   \tilde{C} \tilde{\lambda}^m \inf_{j=0,\ldots,m-1} \left| D\ft_{l,\theta} (\ft^j_{l,\theta} (x)) \right|.
\end{eqnarray*}
Moreover, if $\ft^j_{l,\theta} (x) \not \in W$ for $0\le j \le m$,
then
\begin{eqnarray*}
\left| D\ft^m_{l,\theta} (x) \right| & \ge & \tilde{C} \tilde{\lambda}^m.
\end{eqnarray*}
\end{prop}

We will make use of Koebe's principle, which we quote here.
See \cite{melstr93} for the proof and additional information.
Let $U\subset V$ be two intervals.
We say that $V$ contains a $\delta$-scaled neighborhood of $U$
if both components of $V \backslash U$ have at least length $\delta |U|$.
\begin{defn}
The distortion of a diffeomorphism $g$ on an interval $I$
is defined as
$$
\sup_{x,y \in I} \frac{|Dg(x)|}{|Dg(y)|}.
$$
\end{defn}

\begin{theorem}{\rm \bf [Koebe principle]}\label{koebe}
Let $f$ have negative Schwarzian derivative and let $\Uc$
be a neighborhood of $f$ in $C^3$.
Then for each $\delta >0$, there exists $K < \infty$ so that
$g \in \Uc$, if $I\subset J$ are intervals, 
$\Res{g^n}{J}$ is a diffeomorphism and $g^n(J)$ contains 
a $\delta$-scaled neighborhood of $g^n(I)$ then the distortion of
$\Res{g^n}{I}$ is bounded by $K$. 
\end{theorem}

A useful property of maps with negative Schwarzian derivative 
is the following principle, see \cite{melstr93}.

\begin{theorem}{\rm \bf [Minimum principle]}\label{minimum}
Let $f$ be a map with negative Schwarzian derivative on a closed
interval $I = [a,b]$. If $Df$ does not vanish on $I$, then
\begin{eqnarray*}
|Df(x)| & \ge & \min\{ |Df(a)| , |Df(b)| \}
\end{eqnarray*}
for all $x \in I$.
\end{theorem}

The following lemma is similar to Theorem III.6.2 in \cite{melstr93}.
So is its proof.

\begin{lem}\label{lemma_size}
There are constants $K>0$, $0<\rho<1$, so that for all large enough $l$
the following holds.
Let $I_m$ be a maximal interval
with $\Res{\ft^m_{l,\theta_l^*}}{I_m}$ a homeomorphism.
Then
\begin{eqnarray*}
\left| I_m \right| & \le & K \rho^m.
\end{eqnarray*}
\end{lem}

\noindent {\sc Proof.}
Denote $\tilde{{\cal O}}(c) = \{ \ft_{l,\theta^*_l}^i (c) \}_{i\ge 0}$.
The number of elements in $\tilde{{\cal O}}(c)$ is fixed. There is also 
a minimum distance between any two points in $\tilde{{\cal O}}(c)$, 
uniformly in $l$. Thus we may let $W$ be a neighborhood of $c$ 
such that $\tilde{{\cal O}}(c) \cap W = \emptyset$ for all $l$.
Similarly, $\tilde{{\cal O}}(e) = \{ \ft_{l,\theta^*_l}^i (e) \}_{i\ge 0}$
is a finite set.

Let $J_m$ be a maximal interval on which $\ft_{l,\theta_l^*}^m$ is
a homeomorphism, but $\ft_{l,\theta_l^*}^{m+1}$ not.
Let $\{J_{m+1}^i\}$  be the subintervals of $J_m$ on which 
$\ft_{l,\theta_l^*}^{m+1}$ is a homeomorphism.
The boundary points of $\ft_{l,\theta_l^*}^m (J_m)$
are contained in $\tilde{{\cal O}}(c) \cup \tilde{{\cal O}}(e)$.
Since $\tilde{{\cal O}}(c)$  and $\tilde{{\cal O}}(e)$  finite,
all intervals $\ft_{l,\theta_l^*}^m (J_m),
\ft_{l,\theta_l^*}^{m+1}(J_{m+1}^i)$ have lengths which are bounded
below uniformly in $l$.
Applying the Koebe principle (to $\ft_{l,\theta_l^*}$) one checks that
there is a constant $\tau < 1$ with
\[
 \frac{ \left| J_{m+1}^i \right|}{ \left| J_{m} \right|}
\le \tau.
\]
Further, since $\tilde{{\cal O}}(c)$ and $\tilde{{\cal O}}(e)$ 
are each finite it is clear that 
$\Res{\ft_{l,\theta_l^*}^{m+k} }{J_{m+1}^i}$ is not a 
homeomorphism for some uniformly bounded $k$.
The result follows and it is clear that the constants 
can be chosen uniformly in $l$.
\qed \\

The next proposition discusses expansion properties
of $\ft_{l,\theta_l^*}$, for $\theta = \theta_l^*$.

\begin{prop}\label{prop_ga}
For any small enough neighborhood 
$W$ of $c$,
there are constants $\tilde{C}>0$ and $\tilde{\lambda} >1$, so that
the following holds for all $l$ sufficiently large.
If $\ft^j_{l,\theta^*_l} (x) \not \in W$ for $0\le j \le m-1$, then
\begin{eqnarray*}
\left| D\ft^m_{l,\theta_l^*} (x) \right| & \ge & \tilde{C} \tilde{\lambda}^m.
\end{eqnarray*}
If $\ft^m_{l,\theta_l^*} (x) \in W$, then 
\begin{eqnarray*}
\left| D\ft^m_{l,\theta_l^*} (x) \right| & \ge & \tilde{C} \tilde{\lambda}^m.
\end{eqnarray*}
Without any condition,
\begin{eqnarray*}
\left| D\ft^m_{l,\theta_l^*} (x) \right| & \ge & 
   \tilde{C} \tilde{\lambda}^m \inf_{j=0,\ldots,m-1} 
    \left| D\ft_{l,\theta_l^*} (\ft^j_{l,\theta_l^*} (x)) 
  \right|.
\end{eqnarray*}
\end{prop}

\noindent {\sc Proof.}
Let $W$ be a neighborhood of $c$, small enough  so that
$\ft^{i}_{l,\theta^*_l}(c) \cap W = \emptyset$ for $i >0$.

We claim that there are $\tilde{C}>0$, $\tilde{\lambda} >1$ so that for 
all large enough integers $l$, if $\ft^j_{l,\theta}(x) \not\in W$ 
for $0\le j\le m-1$, then 
$|D\ft^m_{l,\theta_l^*} (x)| \ge \tilde{C}\tilde{\lambda}^m$.
It suffices to show that there exists $M$ with 
$|D\ft_{l,\theta_l^*}^M(x)| > 1$ (compare the proof of Theorem III.3.3 
in \cite{melstr93}). Assume there exist points $x_m \in [0,1]$ with 
$\ft_{l,\theta_l^*}^j (x_m)\not \in W$, $0\le j \le m-1$, and 
$|D\ft_{l,\theta_l^*}^m (x_m)| \le 1$. Since $\ft_{l,\theta_l^*}$ has 
negative Schwarzian derivative we have that on one side of $x_m$, 
$|D\ft_{l,\theta_l^*}^m(x)| \le 1$. Let $H_m$ be the maximal interval 
bounded by $x_m$ so that $\breve{f}_{l,\theta_l^*}^m(\cdot;H_m)$ is a 
homeomorphism on $H_m$ and 
$\left| D\breve{f}_{l,\theta_l^*}^m (\cdot;H_m) \right| \le 1$ on $H_m$.
This implies that $|  \breve{f}^{m}_{l,\theta_l^*} (H_m;H_m) | \leq |H_m|$.
Let $y_m$ be the other boundary point of $H_m$.
Then either $\breve{f}^{k(m)}_{l,\theta_l^*} (y_m;H_m) = c$
or $\breve{f}^{k(m)}_{l,\theta_l^*} (H_m;H_m) = e$, for some $k(m)$.
In both cases, $| \breve{f}^{k(m)}_{l,\theta_l^*} (H_m;H_m) |$ is bounded away from
$0$.
By lemma~\ref{lemma_size}, $|H_m| \le \tilde{C} \rho^m$ 
is small as $m$ is big.
Since  $|  \breve{f}^{m}_{l,\theta_l^*} (H_m;H_m) | \leq |H_m|$,
also  $\breve{f}^{m}_{l,\theta_l^*} (H_m;H_m) $ is small for large $m$.
The interval $ \breve{f}^{k(m)}_{l,\theta_l^*} (H_m;H_m)$ has length
bounded away from
$0$, say $ |  \breve{f}^{k(m)}_{l,\theta_l^*} (H_m;H_m)| \ge  \delta$.
It follows that $m - k(m)$ tends to $\infty$ as $m \to \infty$.
In fact, if $D$ denotes the minimum of $|D\breve{f}_{l,\theta_l^*}(\cdot;\cdot)|$ 
over $[0,1] \backslash W$, then 
\begin{eqnarray}\label{largerbound}
m  - k(m) & \ge & \ln (\tilde{C} \rho^m / \delta) / \ln D.
\end{eqnarray}
Lemma~\ref{lemma_size} yields that
$|  \breve{f}^{k(m)}_{l,\theta_l^*} (H_m;H_m)| \leq \tilde{C} \rho^{m-k(m)}$,
so that
\begin{eqnarray}\label{smallerbound}
m - k(m) & \le & \ln (\delta/\tilde{C}) / \ln \rho.
\end{eqnarray}
Now (\ref{largerbound}) and (\ref{smallerbound}) contradict each other,
proving the claim.
Observe that $M$ and hence 
the constants $\tilde{C}$ and $\tilde{\lambda}$ can be chosen uniformly in $l$.
This proves the first estimate.

Next, suppose $\ft_{l,\theta^*_l} (x) \in W$.
Let $H_m$ be the maximal interval containing $x$ such that
$\breve{f}^m_{l,\theta} (\cdot; H_m)$ is a homeomorphism
on $H_m$.
Because the orbit of $c$  is finite, the interval
$\breve{f}^m_{l,\theta} (H_m; H_m)$ extends a positive distance away from
$\breve{f}^m_{l,\theta} (x; H_m)$ to both sides.
Koebe's principle implies 
$| D\breve{f}_{\eta_l}^j (x;H_m)| \ge \tilde{C} \tilde{\lambda}^j$ for some $\tilde{C}>0, \tilde{\lambda} >1$
which gives the second estimate.

Finally, not assuming any condition, split the iterates $x, \ldots, \ft_{l,\theta}^m(x)$ 
into a part that ends in $W$, one iterate starting in $W$, and a part that stays outside $W$.
Combining the first two estimates for the first and last part, proves the last estimate.
\qed \\

\noindent {\sc Proof of proposition~\ref{prop_exp}.}
As in the proof of Theorem III.6.4 in \cite{melstr93}.
\qed\\

\subsection{Parameter dependence}\label{sec_par}

\begin{prop}\label{prop_parameterstate}
There are constants $C,\varepsilon,K >0$ so that for all large enough $l$
the following holds.
If $|\theta - \theta_l^*| < \varepsilon$, $k > K$ and
\begin{eqnarray*}
| (\ft^k_{l,\theta})' (f_{l,\theta}(c)) | & \geq & e^{\tilde{\lambda} k},
\end{eqnarray*}
then for all $k>K$,
\[
\frac{1}{C} 
\le 
\frac{ | \frac{d}{d \theta} \ft^k_{l,\theta}(c) |}
     { | (\ft^{k-1}_{l,\theta})' (f_{l,\theta}(c)) | }
\le
C.
\]
\end{prop}

\noindent {\sc Proof.}
Let $x(l,\theta)$ be the continuation of $\ft_{l,\theta_l^*}(c)$ 
for $\theta$ near $\theta_l^*$
given by $\ft^{k - 1}_{{l,\theta}} (x(l,\theta)) = y^*_{l,\theta}$
(where $y^*_{l,\theta^*}$ is the hyperbolic periodic point from the definition
of $\theta_l(\ga_l^*) = \theta_l^*$).
Then Proposition~\ref{localmap} implies that
\begin{eqnarray}\label{estder}
\frac{1}{K} \le 
\left|\frac{d}{d \theta} (x(l,\theta) - \ft_{l,\theta} (c))
\right| 
\le K,
\end{eqnarray}
for some constant $K>1$.
Writing
$\ft_{l,\theta}^{k-1} (\ft_{l,\theta}(c)) = 
\ft_{l,\theta}^{k-1} (x(l,\theta)
+ \ft_{l,\theta}(c) - x(l,\theta))$,
the chain rule gives
\begin{eqnarray*}
\Res{\frac{d}{d \theta} \ft_{l,\theta}^{k-1} 
   (\ft_{l,\theta}(c))}{\theta = \theta_l^*}
& = &
\Res{\frac{d}{d \theta} \ft_{l,\theta}^{k-1}(x)}
            {\theta = \theta_l^*,x = x(l,\theta)}
+
(\ft_{l,\theta_l^*}^{k-1})' (\ft_{l,\theta_l^*}(c))
\Res{\frac{d}{d \theta} ( \ft_{l,\theta}(c) - x(l,\theta))}{\theta = \theta_l^*}.
\end{eqnarray*}
Now 
$\Res{\frac{d}{d \theta} \ft_{l,\theta}^{k-1} (x_l(l,\theta))}{\theta = \theta_l^*} 
\equiv \frac{d}{d \theta} y^*_{l,\theta}$
is arbitrarily small for $l$ sufficiently large.
By (\ref{estder}) and the exponential growth
of $(\ft_{l,\theta_l^*}^{k-1})' (\ft_{l,\theta_l^*}(c))$,
the statement of the proposition holds for $\theta = \theta_l^*$.

By the chain rule,
\begin{eqnarray*}
(\ft^{k-1}_{l,\theta})' (\ft_{l,\theta}(c)) & = &
\ft_{l,\theta}' (\ft_{l,\theta}^{k-1}(c)) 
  (\ft_{l,\theta}^{k-2})' (\ft_{l,\theta}(c)),
\\
\frac{d}{d \theta} \ft_{l,\theta}^k(c) & = &
\ft_{l,\theta}' (\ft_{l,\theta}^{k-1}(c)) 
 \frac{d}{d \theta} \ft_{l,\theta}^{k-1} (c)
 + 
\frac{\partial}{\partial \theta} \ft_{l,\theta}(\ft_{l,\theta}^{k-1} (c)).
\end{eqnarray*}
It follows that
\begin{eqnarray*}
\left|
\frac{\frac{d}{d \theta} \ft_{l,\theta}^k(c)}
     {(\ft^{k-1}_{l,\theta})' (\ft_{l,\theta}(c))}
 -
\frac{\frac{d}{d \theta} \ft_{l,\theta}^{k-1}(c)}
     {(\ft^{k-2}_{l,\theta})' (\ft_{l,\theta}(c))}
\right|
& = &
\left|
\frac{
\ft_{l,\theta}' (\ft_{l,\theta}^{k-1}(c)) 
 \frac{d}{d \theta} \ft_{l,\theta}^{k-1} (c)
 + 
\frac{\partial}{\partial \theta} \ft_{l,\theta}(\ft_{l,\theta}^{k-1} (c))
}
{
\ft_{l,\theta}' (\ft_{l,\theta}^{k-1}(c)) 
  (\ft_{l,\theta}^{k-2})' (\ft_{l,\theta}(c))}
-
\frac{\frac{d}{d \theta} \ft_{l,\theta}^{k-1}(c)}
     {(\ft^{k-2}_{l,\theta})' (\ft_{l,\theta}(c))}
\right|
\\
& = &
\left|
\frac{\frac{\partial}{\partial \theta} \ft_{l,\theta}(\ft_{l,\theta}^{k-1} (c))}
     {(\ft_{l,\theta}^{k-1})' (\ft_{l,\theta}(c))}
\right|.
\end{eqnarray*}
By assumption,
$| {(\ft_{l,\theta}^{k})' (\ft_{l,\theta}(c))}| \ge e^{\tilde{\lambda} k}$ for $k>K$.
For each positive integer $K$, there are  constants $C, \varepsilon$ so that
\[
\frac{1}{C}
\le
\left|
\frac{ \Res{\frac{d}{d \theta} \ft_{l,\theta}^{K}(c)}{\theta = \theta_l^*}}
     {(\ft^{K-1}_{l,\theta_l^*})' (\ft_{l,\theta_l^*}(c))}
\right|
 \le  C,
\]
if $| \theta - \theta_l^* | \le \varepsilon$.
Hence,
\begin{eqnarray*}
\left|
\frac{\frac{d}{d \theta} \ft_{l,\theta}^k(c)}
     {(\ft^{k-1}_{l,\theta})' (\ft_{l,\theta}(c))}
 -
\frac{\frac{d}{d \theta} \ft_{l,\theta}^{K}(c)}
     {(\ft^{K-1}_{l,\theta})' (\ft_{l,\theta}(c))}
\right|
& \le &
C_K,
\end{eqnarray*}
for  $| \theta - \theta_l^* | \le \varepsilon$.
The constant $C_K$ is small if $K$ is large.
The proposition follows.
\qed\\

\subsection{Binding}\label{sec_bin}

The next proposition is also used in Section~\ref{sec_measure}
to show the existence of absolutely continuous invariant measures
of $f_{l,\theta}$ for $\theta \in \Theta$.
The binding period of $\omega \in \Pc_l^{(n)}$ is defined by (\ref{binding1}), (\ref{binding2}).

\begin{prop}\label{lem_binding}
There exist constants $\tilde{\la}$ and $\beta<1$ such 
that the following holds.
Let $0\le k \le n-1$, $\omega \in \Pc_l^{(n)}$
and suppose that $k$ is an essential return time for $\omega$
with return depth $r$. Then the binding period
$p = p(\omega,k)$ satisfies 
 $p \le 2r$. 
We have
\begin{eqnarray*}
\left| D\ft^{p+1}_{l,\theta} (\ct_k({l,\theta})) \right| & \ge & 
e^{\frac{1}{6} \tilde{\lambda} (p+1)}. 
\end{eqnarray*}
Furthermore,
\begin{eqnarray*}
|\omega_{k+p+1}| & \ge & |\omega_k|^\beta.
\end{eqnarray*}
\end{prop}

Recall that $\eta_0 = (c , \ct_k({l,\theta}))$, or 
$(\ct_k({l,\theta}),c)$, and
$\eta_j =  \breve{f}_{l,\theta}^j (\eta_0;\eta_0)$.
The next lemma implies bounded distortion of iterates of $\breve{f}_\ga$
on $\eta_0$ during the binding period. The lemma is an ingredient for the proof of
Proposition~\ref{lem_binding}, as in \cite{luz00}. The proof of the lemma differs
from \cite{luz00}; distortion estimates for a passage through $\Et_\ga$
have to be treated separately. 
The remainder of the proof of Proposition~\ref{lem_binding}
is as in \cite{luz00}.

\begin{lem}
For $y_0,z_0 \in \eta_0$,
\begin{eqnarray*}
\left|
\frac{ D \breve{f}_{l,\theta}^i (y_0;\eta_0)}{D \breve{f}_{l,\theta}^i (z_0;\eta_0)}
\right|
& \le & K,
\end{eqnarray*}
for $0\le i \le p$ and 
for some $K>0$.
\end{lem}

\noindent {\sc Proof.}
Write
$y_j = \breve{f}_{l,\theta}^j (y_0;\eta_0)$ and
$z_j = \breve{f}_{l,\theta}^j (z_0;\eta_0)$.
By the chain rule,
\begin{eqnarray*}
\left|
\frac{ D \breve{f}_{l,\theta}^i (y_0;\eta_0)}{D \breve{f}_{l,\theta}^i (z_0;\eta_0)}
\right|
& = &
\left|
\prod_{j=0}^{i-1} 
\frac{ D \breve{f}_{l,\theta} (y_j;\eta_0)}{D \breve{f}_{l,\theta} (z_j;\eta_0)}
\right|
\\
& = &
\left|
\prod_{j=0}^{i-1} 
\left(
1 + \frac{ D \breve{f}_{l,\theta} (y_j;\eta_0) - D \breve{f}_{l,\theta} (z_j;\eta_0)}
         {D \breve{f}_{l,\theta} (z_j;\eta_0)}
\right) \right|.
\end{eqnarray*}
If $\eta_j$ lies outside $\Et_{l,\theta}$, then
$| D^2 \breve{f}_{l,\theta} (\cdot;\eta_0) |$ is bounded by a constant $C>0$
and
$$ | D \breve{f}_{l,\theta} (y_j;\eta_0) - D \breve{f}_{l,\theta} (z_j;\eta_0)|
    \le C |\eta_j|.
$$
If $z_j, y_j \in \Et_{l,\theta}$, 
we can write $z_j = \phi^s_{l,\theta} (y_j)$.
Thus $\breve{f}_{l,\theta} (z_j;\eta_0)= 
\phi^{-s}_{l,\theta} \circ \breve{f}_{l,\theta} (  \phi^s_{l,\theta} (y_j) ;\eta_0)$.
Using this, it follows that again
$ | D \breve{f}_{l,\theta} (y_j;\eta_0) - D \breve{f}_{l,\theta} (z_j;\eta_0)|
\le C |\eta_j|$
for some $C>0$.
Hence,
\begin{eqnarray*}
\left|
\frac{ D \breve{f}_{l,\theta}^i (y_0;\eta_0)}{D \breve{f}_{l,\theta}^i (z_0;\eta_0)}
\right|
& \le &
e^{\ln \prod_{j=0}^{i-1} 
\left( 
1 + C \left| \frac{\eta_j}{ D \breve{f}_{l,\theta} (z_j;\eta_0)} \right|
\right)}
\\
& = &
e^{ \sum_{j=0}^{i-1}
\ln \left( 
1 + C \left| \frac{\eta_j}{ D \breve{f}_{l,\theta} (z_j;\eta_0)} \right|
\right)}
\\
& \le &
e^{ C \sum_{j=0}^{i-1} \left| \frac{\eta_j}{ D \breve{f}_{l,\theta} (z_j;\eta_0)} 
      \right|}.
\end{eqnarray*}
We proceed to estimate $ 
\sum_{j=0}^{i-1} \left| \frac{\eta_j}{ D \breve{f}_{l,\theta} (z_j;\eta_0)} 
      \right|$.
By the definition of binding period, 
\begin{equation}\label{etaj}
|\eta_j| \le e^{-2\alpha j}.
\end{equation}
The bounded recurrence assumption implies that 
$|\tilde{c}_k - c| \ge e^{-\alpha k}$.
Write $\breve{f}_{l,\theta}^j (c;\eta_0) =   \tilde{c}_{k(j)}$.
If $\Et_{l,\theta}$ is small enough, we have $k(j) \le \frac{4}{3} j$ 
(see Section~\ref{sec_relate}).
We may hence assume that
$|\breve{f}_{l,\theta}^j (c;\eta_0) - c| \ge e^{-\frac{4}{3} \alpha j}$.
It follows from this and (\ref{etaj})  that
$|z_j - c| \ge e^{-\frac{4}{3} \alpha j} - e^{-2 \alpha j}
=  e^{-\frac{4}{3} \alpha j} ( 1 - e^{-\frac{2}{3} \alpha j})$.
Further $|D \breve{f}_{l,\theta} (z_j;\eta_0) | \ge C |z_j - c|$ for some $C>0$,
so that
\begin{equation}\label{zj}
|D \breve{f}_{l,\theta} (z_j;\eta_0) | \ge
C e^{-\frac{4}{3} \alpha j},
\end{equation}
for some $C>0$.
Combining (\ref{etaj}) and (\ref{zj}) shows that
$ \sum_{j=0}^{i-1} \left| \frac{\eta_j}{ D \breve{f}_{l,\theta} (z_j;\eta_0)} 
   \right|$
is bounded, thus proving the lemma.
We remark that the distortion bound is close to 1 
if $|\theta - \theta^*_l|$ is small. This follows from the observation that 
$\breve{f}_{l,\theta}^j (\eta_0; \eta_0)$  is outside a neighborhood of $c$, 
where it undergoes
exponential expansion, for a large number of iterates.
\qed

\subsection{Induction}

The main results from the previous sections are
Propositions~\ref{prop_exp}, \ref{prop_parameterstate}, \ref{lem_binding}.
These results have their counterparts in proofs of
the work of Benedicks and Carleson. From this point on, we can
follow \cite{luz00} closely. For completeness we sketch the 
remaining steps leading to the proof of Proposition~\ref{jakobson}
in the next two sections. 

In the inductive constructions, the following two propositions 
are shown to hold.
The proofs are as in \cite{luz00}, relying on 
Propositions~\ref{prop_exp}, \ref{prop_parameterstate}, and \ref{lem_binding}.

\begin{prop}{\bf [Bounded recurrence]}
Each point in $\Theta_l^{(n)}$ satisfies $(\BR)_n$.
\end{prop}

One shows that in fact 
${\cal R}_l^{(k)}(\theta) \le 5 {\cal E}_l^{(k)} (\theta)$, 
$0\le k \le n$, if $\theta \in \Theta_l^{(n)}$. 
That is, a substantial proportion of the returns
are chopping times. By assumption, ${\cal E}^{(k)} (\ga) \le \alpha k /10$.
Hence, ${\cal R}^{(k)} (\ga) \le \alpha k /2$.
As in \cite{luz00} one shows that this bound implies $(\BR)_n$.

\begin{prop}{\bf [Bounded distortion]}
Restricted to a connected component of $\omega \in \Pc_l^{(n)}$, 
the map $\ct_j$ is a diffeomorphism
with uniformly bounded distortion for all $j\le \nu + p + 1$ where $\nu \le n$
is the last essential or inessential return time of $\omega$ and $p$
is the associated binding period.
If $n > \nu+p+1$ then the same statement holds for all $j\le n$ restricted
to any subinterval $\bar{\omega}$ such that $\bar{\omega_j} \subset \Delta^+$. 
\end{prop}

\subsection{Combinatorial estimates and measure bounds}

For the proof of Proposition~\ref{jakobson}, combinatorial 
properties of $\Pc_l^{(n)}$ are studied.
Escape times play a central role in this study.
The combinatorial properties described next are used at the end of the
section to prove Proposition~\ref{jakobson}.

To each $\omega \in \hat\Pc_l^{(n)}$ is associated a sequence
$0=\eta_0 < \eta_1 < \ldots < \eta_s \le n$, $s = s(\omega)$ of escape times
and a corresponding sequence of escaping components
$\omega \subset \omega^{(\eta_s)} \subset \ldots \subset \omega^{(\eta_0)}$ 
with $\omega^{(\eta_i)} \in \Pc_l^{(\eta_i)}$.
Let $\omega^*_i = \omega^{(\eta_i)}$ for $1 \le i \le s$ and 
$\omega^*_i = \omega$ for $s+1 \le i \le n$.
This defines $\omega^*_i$ for each $0\le i \le n$.
Observe that for $\omega,\nu \in \Pc_l^{(n)}$
and $0\le i \le n$,
the sets $\omega^*_i$ and $\nu^*_i$ are either disjoint
or coincide.
Define
\begin{eqnarray*}
Q_l^{(i)} & = & \bigcup_{\omega \in \hat{\Pc_l}^{(n)}} \omega^*_i
\end{eqnarray*}
and let 
\begin{eqnarray*}
\Qc_l^{(i)} = \{ \omega^*_i \}
\end{eqnarray*}
be the natural partition of $Q_l^{(i)}$ into sets of the form $\omega^*_i$.
Observe that $\Theta_l^{(n-1)} = Q_l^{(n)} \subset \ldots \subset Q_l^{(0)} = \Theta_l^{(0)}$
and $\Qc_l^{(n)} = \hat\Pc_l^{(n)}$.
For $\omega = \omega^*_i \in \Qc_l^{(i)}$, $0\le i \le n-1$, let
\begin{eqnarray*}
Q_l^{(i+1)} (\omega) & = & \left\{ \omega' = \omega^*_{i+1} \in \Qc_l^{(i+1)}; \;\;
                             \omega' \subset \omega \right\}.
\end{eqnarray*}
Denote by $\Qc_l^{(i+1)} (\omega)$ the partition $\Qc_l^{(i+1)}$ 
restricted to $\omega$.
Define $\Dc_l^{(i+1)} : Q_l^{(i+1)} (\omega) \to \bN$ by 
\begin{eqnarray*}
\Dc_l^{(i+1)} (\omega') & = & {\cal E}_l^{(\eta_{i+1})}(\omega') - 
   {\cal E}_l^{(\eta_{i})}(\omega')
\quad \text{for} 0\le i \le s,
\end{eqnarray*}
and $\Dc_l^{(i+1)} =0$ for $i>s$. 
Let
\begin{eqnarray*}
Q_l^{(i+1)} (\omega,R) & = & \left\{
                \omega' \in Q_l^{(i+1)}; \;\; \omega' \subset \omega,
      \Res{  \Dc_l^{(i)} }{\omega'} = R 
                          \right\}.
\end{eqnarray*}

\begin{prop}\label{comb}
\begin{eqnarray*}
\sum_{\omega' \in \Qc_l^{(i+1)}(\omega,R)}  |\omega'| 
& \le &
e^{ - \eta R} |\omega|,
\end{eqnarray*}
for some $\eta > 0$.
\end{prop}

\noindent {\sc Proof.}
One can take $- \eta = 10 \frac{\alpha}{\tilde{\lambda}} - 1$, which is negative
for $\alpha$ small enough.
The proof divides into two parts.
One bounds the cardinality of $\Qc_l^{(i)} (\omega,R)$ by $e^{\beta R}$
and one shows that for any $\omega \in \Qc_l^{(i)}$, $0\le i \le n-1$, and
$\hat{\omega} \in \Qc_l^{(i+1)} (\omega, R)$, one has
$|\hat{\omega} | \le e^{ (9 \beta - 1) R} | \omega|$.
Combining the two statements proves the proposition.
For the proofs one can follow \cite{luz00}. 
\qed \\

\begin{lem}\label{lem_estcalE}
$$
\int_{\Theta_l^{(n-1)}} e^{{\cal E}^{(n)}/2}
    = \sum_{\omega \in \Qc_l^{(n)}}  e^{{\cal E}^{(n)}/2} |\omega|
    \le e^{3n/r_\delta} |\Theta_l^{(0)}|.
$$
\end{lem}

\noindent {\sc Proof.}
The equality follows immediately from the definitions. 
For the
inequality, let $0\le i \le n$, $\omega \in \Qc_l^{(i)}$ and write
\begin{eqnarray*}
\sum_{\omega' \in \Qc_l^{(i+1)}(\omega)} e^{\Dc_l^{(i)} (\omega')/2} |\omega'|
& = & 
\sum_{\omega' \in \Qc_l^{(i+1)} (\omega,0)} |\omega'| + 
\sum_{R \ge r_\delta} e^{R/2} \sum_{\omega' \in \Qc_l^{(i+1)}(\omega,R)} |\omega'|.
\end{eqnarray*}
Proposition~\ref{comb} (with the remark from its proof that $-\eta = 10 \beta - 1$ with
$\beta = \alpha/\tilde{\lambda}$) implies
\begin{eqnarray*}
\sum_{\omega' \in \Qc_l^{(i+1)}(\omega)} e^{\Dc_l^{(i)} (\omega')/2} |\omega'|
+ \sum_{R \ge r_\delta} e^{R/2} \sum_{\omega' \in \Qc_l^{(i+1)}(\omega,R)} |\omega'|
& \le &
\left( 1 + \sum_{R\ge r_\delta} e^{(10 \beta - \frac{1}{2}) R}\right) |\omega|,
\end{eqnarray*}
so that
\begin{eqnarray}\nonumber
\sum_{\omega' \in \Qc_l^{(i+1)}(\omega)} e^{\Dc_l^{(i)} (\omega')/2}|\omega'|
& \le &
(1 + e^{-r_\delta/3}) |\omega|
\\
\label{repeatedly}
& \le & e^{3/r_\delta} |\omega|,
\end{eqnarray}
assuming that $\beta$ has been chosen small enough and $r_\delta$ large enough.
Since
${\cal E}^{(n)} = \Dc_l^{(0)} + \ldots + \Dc_l^{(n-1)}$
and $\Dc_l^{(i)}$ is constant on elements of $\Qc_l^{(i)}$, we have
\begin{eqnarray*}
\lefteqn{ \sum_{\omega \in \Qc_l^{(n)}} e^{{\cal E}^{(n)} (\omega)/2} |\omega| = }
\\
&  &
\sum_{\omega^*_1} \in \Qc_l^{(1)}(\omega^*_2) e^{\Dc_l^{(0)} (\omega^*_1)/2}
\cdots
\sum_{\omega^*_{n-1} \in \Qc_l^{(n-1)}(\omega^*_n)} e^{\Dc_l^{(n-2)} (\omega^*_{n-1})/2}
\sum_{\omega = \omega^*_n \in \Qc_l^{(n)}} e^{\Dc_l^{(n-1)} (\omega^*_n)/2} |\omega|.
\end{eqnarray*}
Applying (\ref{repeatedly}) repeatedly gives
\begin{eqnarray*}
\sum_{\omega \in \Qc_l^{(n)}} e^{{\cal E}^{(n)} (\omega) / 2} |\omega|
& \le &
e^{3 n /r_\delta} |\Theta_l^{(0)}|.
\end{eqnarray*}
\qed\\
     
Observe
\begin{eqnarray*}
\left| \Theta_l^{(n-1)} \backslash \Theta_l^{(n)} \right|
& = &
\left| \{ \omega \in \Qc_l^{(n)}; \;\; e^{{\cal E}^{(n)} /2} \ge
    e^{\alpha n / 20} \} \right|.
\end{eqnarray*}
Chebyshev's inequality and Lemma~\ref{lem_estcalE} yield
\begin{eqnarray*}
\left| \Theta_l^{(n-1)} \backslash \Theta_l^{(n)} \right|
& \le &
e^{- \alpha n / 20} \int_{\Theta_l^{(n-1)}} e^{{\cal E}^{(n)} /2}
\\
& \le &
e^{\left( \frac{3}{r_\delta} - \frac{\alpha}{20} \right) n} |\Theta_l^{(0)}|
\\
& \le & e^{-\alpha n /30} |\Theta_l^{(0)}|,
\end{eqnarray*}
if $r_\delta$ is large enough.
This implies
\begin{eqnarray*}
|\Theta_l^{(n)} | & \ge & |\Theta_l^{(n-1)} | - e^{-\alpha n /30}
 |\Theta_l^{(0)}|.
\end{eqnarray*}
Write $\Theta_l^{(0)} = (\theta_l^*  - \epsilon, \theta_l^* + \epsilon)$.
For $\epsilon$ small, 
there exists $N$ so that
$\Theta_l^{(j)} = \Theta_l^{(j+1)}$ for all $j \le N$.
Hence
\begin{eqnarray*}
|\Theta_l^{(n)}| & \ge & \left( 1 - \sum_{i=N}^n e^{- \alpha i/30}
 \right) |\Theta_l^{(0)}|.
\end{eqnarray*}
Noting that $|\Theta_l^{(0)}| = 2 \epsilon$, for all $l$, 
it follows that a uniform lower bound for $m(\Theta_l)$ exists.
This concludes the proof of Proposition~\ref{jakobson}.
Since $N$ goes to $\infty$ as $\epsilon \to 0$,
we have also shown that $\theta^*_l$ is a Lebesgue density 
point of $\Theta_l$.

\subsection{Proof of Theorem~\ref{existence}}\label{sec_proof}

In the previous sections it was proved that 
$\ft_{l,\theta}$ has bounded recurrence for a set of parameter 
values $\theta \in \Theta_l$ with measure bounded from below, 
uniformly in $l$. By Proposition~\ref{density}, 
this implies that $\Omega = \cup_l g_l(\Theta_l)$ has positive measure
and has positive density at $\ga = 0$.

The following proposition implies that $f_{l,\theta}$, 
$\theta \in \Theta_l$, has exponential expansion along the 
orbit of $f_{l,\theta} (c)$.
The proof is as in \cite{luz00}.

\begin{prop}\label{prop_br>exp}
If $\theta$ is close to $\theta_l^*$ and $\ft_{l,\theta}$ satisfies 
$(\BR)_k$, then
\begin{eqnarray*}
\left| D\ft^{k+1}_{l,\theta} (f_{l,\theta} (c)) \right| & \ge & 
\tilde{C} e^{\tilde{\lambda} (k+1)},
\end{eqnarray*}
for some $\tilde{C} > 0$, $\tilde{\lambda} > 1$.
\end{prop}

Combining Proposition~\ref{prop_br>exp} with 
Lemma~\ref{lemma_withouttilde},
gives

\begin{prop}\label{prop_coleck}
For each $\ga \in \Omega$, there are $\tilde{C}>0$, $\lambda > 1$, so that
\begin{eqnarray*}
\left| Df^{k}_\ga (f_\ga (c)) \right| & \ge & 
\tilde{C} e^{\lambda k}.
\end{eqnarray*}
\end{prop}

Thus $f_\ga$ is a Collet-Eckmann map if $\ga \in \Omega$.
Collet-Eckmann maps are known to admit
absolutely continuous invariant measures, see Theorem V.4.6 in \cite{melstr93}.
This concludes the proof of Theorem~\ref{existence}, except for the
conclusion that $\mathrm{supp}(\nu_\ga) = [f^2_\ga(c),f_\ga(c)]$,
which we postpone until the next section.


\section{Intermittency}\label{sec_measure}

In this section we study intermittent time series
of $\ft_\ga$ 
at parameter values $\ga$ for which $f_\ga$ admits an
absolutely continuous invariant measure.
A proof of Theorem~\ref{intermittency} is in Section~\ref{sec_inter}.
For the proof one needs to know the
(average) length of the relaminarization phase.
The relaminarization phase is studied in Section~\ref{app}.

\subsection{Invariant measures and intermittency}\label{sec_inter}

In the previous section we constructed a set $\Omega$ of parameter values
with positive density at $\ga = 0$, so that $f_\ga$ has bounded recurrence
(see Definition~\ref{def_br}) for $\ga \in \Omega$.
By Proposition~\ref{prop_coleck}, the bounded recurrence condition implies
the Collet-Eckmann condition stating that $|Df^n_\ga (f_\ga(c))| \ge C \la^n$ 
for some $C>0$, $\la > 1$.
It is known that a map $f_\ga$ satisfying the Collet-Eckmann condition,
possesses an absolutely continuous invariant measure, see \cite{melstr93}.
Because we need bounds on the density of the invariant measures 
in our discussion of intermittency, we give an alternative way to
produce invariant measures following \cite{you92}. (See also \cite{ryc88} and
\cite{rycsor92})

\begin{prop}\label{prop_tildeacip}
For $\ga \in \Omega$, $\ft_\ga$ possesses an absolutely
continuous invariant measure $\tilde{\nu}_\ga$.
There is a constant $K>0$ not depending on $\ga$ so that for
any Borel set $A$, 
\begin{eqnarray*}
\tilde{\nu}_\ga (A) & \le & K \sqrt{m(A)}.
\end{eqnarray*}
The support of $\tilde{\nu}_\ga$ equals $[f^2_\ga(c) , f_\ga(c)]$.
\end{prop}

\noindent {\sc Proof.}
By Proposition~\ref{lem_binding}, the binding period $p(x)$ defined near $c$
satisfies
\begin{description}
\item[(i)] $p(x)  \le  - \tilde{C} \ln |x-c|$ for some $\tilde{C}>0$ not 
depending on $\ga$,
\item[(ii)] $|D\ft_\ga^j (f_\ga(x))| \ge \tilde{C} \tilde{\lambda}^j$ for some $\tilde{C}>0$,
$\tilde{\lambda} > 1$ and for all $0\le j < p$,
\item[(iii)] $|D\ft_\ga^{p} (x)| \ge \sigma^p$ 
for some $\sigma >1$.
\end{description}
Outside the domain of definition of $p$, let $p = 1$.
Define the return map $R_\ga$ on $[0,1]$ by
\begin{eqnarray*}
R_\ga(x) & = & \ft_\ga^{p(x)}
\end{eqnarray*}
For any $K>0$ and all $\ga \in \Omega$ we may 
assume that $\Res{p}{\Delta}$ is bounded from below by $K$,  
by taking $\Delta$ small enough.
By Proposition~\ref{prop_exp}, if $x,\ldots,\ft_\ga^{m-1}(x) \not\in \Delta$ 
and $\ft^m_\ga(x) \in \Delta$, then $|D\ft^m_\ga (x) | \ge \tilde{C} \tilde{\lambda}^m$
for some $\tilde{C} > 0, \tilde{\lambda} > 1$.
As stated in Proposition~\ref{prop_exp}, 
the constant $\tilde{C}$ does not depend on $\Delta$.
By taking $\Delta$ small, $\tilde{C} \tilde{\la}^K >1$.
It follows that some power of $R_\ga$ is expanding.
The first part of the proposition 
follows as in the proof of Theorem~1 in \cite{you92}.
We sketch the necessary arguments to make clear that the constant $K$ 
does not depend on $\ga$. Details are left to the reader.
The measure $\tilde{\nu}_\ga$ is constructed by finding its density 
as a fixed point of a Perron-Frobenius operator $P_\ga$.
Define
$$
P_\ga \phi (x) = \sum_{y \in \ft_\ga^{-1} (x)} \frac{1}{(\ft_\ga)' (x)} \phi (y).
$$
Note that 
$$
P^n_\ga \phi (x) = \sum_{y \in \ft_\ga^{-n} (x)} \frac{1}{(\ft^n_\ga)' (x)} \phi (y).
$$
Let $g_n : [0,1] \to \bR$ be the function
$g_n (x) = 1 / (\ft^n_\ga)' (x)$ where $\ft^n_\ga$ is continuous and 
$g_n (x) = 0$
elsewhere.
One shows that the variation of $g_n$ is bounded, uniformly in $n$ and $\ga$.
This relies on the negative Schwarzian derivative (see \cite{you92}) and
the analysis of the local saddle-node bifurcation in Section~\ref{sec_snlocal}
(see Proposition~\ref{embedding}). 
By induction the variation of $g_n$ is bounded.
Using the uniform bound on the variation of $g_n$, 
one bounds the variation of $P^n_\ga \phi$ for  
a density $\phi$ with bounded variation.
It follows that $P_\ga$ has a fixed point with uniformly bounded variation.

If $\zeta_\ga$ denotes the measure whose density is the fixed point of $P_\ga$,
then $\tilde{\nu}_\ga$ is obtained by pushing forward $\zeta_\ga$,
$$
\tilde{\nu}_\ga  (A) = 
\sum_{k=1}^\infty  \sum_{j=0}^{k-1} 
\zeta_\ga (\ft^{-j}_\ga (A) \cap B_k),   
$$
where $B_k$ is the set on which $R_\ga = \ft_\ga^k$.
The uniform bound for $\tilde{\nu}_\ga$ follows from the properties of $R_\ga$,
as in \cite{you92}, see also \cite{homyou00}.

That the support of $\tilde{\nu}_\ga$ equals all of $[f^2_\ga(c) , f_\ga(c)]$
follows as in Theorem~2 in \cite{you92}.
It relies on the fact that, if $\ga \in \Omega$, then for each interval $I \subset [0,1]$,
there exists $n$ so that $[f^2_\ga(c),f_\ga(c)] \subset f_\ga^n(I)$.
The necessary topological expansion at $\ga = 0$ is 
guaranteed by Lemma~\ref{lem_over}.
\qed\\

The invariant measure $\nu_\ga$ for $f_\ga$ is constructed
by pushing forward $\tilde{\nu}_\ga$.
Recall that $\Et_\ga = \cup_{i=0}^l
E^{-i}$ with $E^{-i} = ( \Res{f}{E} )^{-iq} (I^u_{\ga})$.
Let
\begin{eqnarray*}
\bar{\nu}_\ga (A) & = & \tilde{\nu}_\ga (A \cap ([0,1]\backslash \Et_\ga)) + 
\sum_{k=0}^{l} \sum_{i=0}^{kq-1} 
       \tilde{\nu}_\ga (f_\ga^{-i}(A) \cap E^{-k}).
\end{eqnarray*}
This measure is obviously finite and can thus be rescaled to a probability 
measure $\nu_\ga$. 

\begin{lem}\label{prop_mb}
The measure $\nu_\ga$, $\ga \in \Omega$, 
is an absolutely continuous invariant probability measure for $f_\ga$. 
\end{lem}

\noindent {\sc Proof.} 
To prove that $\nu_\ga$ is invariant for $f_\ga$,
recall that $\Res{ \ft_\ga}{E^{-j}} = f^{jq}_\ga$ and
$\ft_\ga$ equals $f_\ga$ outside $\Et_\ga$.  
Compute
\begin{eqnarray}
\nonumber
\bar{\nu}_\ga (f^{-1}_\ga(A)) & = & 
\tilde{\nu}_\ga (f^{-1}_\ga(A) \cap  ([0,1]\backslash \Et_\ga)) +
\sum_{k=0}^{l} \sum_{i=1}^{kq-1} \tilde{\nu}_\ga  (f_\ga^{-i}(A) \cap E^{-k})
+
\\
\nonumber
& & 
\sum_{k=0}^{l} \tilde{\nu}_\ga  (f_\ga^{-kq}(A) \cap E^{-k})
\\
\nonumber
& = & 
\tilde{\nu}_\ga (\ft^{-1}_\ga(A) \cap  ([0,1]\backslash \Et_\ga)) +
\sum_{k=0}^{l} \sum_{i=1}^{kq-1} \tilde{\nu}_\ga  (f_\ga^{-i}(A) \cap E^{-k})
+
\\
\label{eqn_inv}
& & 
\sum_{k=0}^{l} \tilde{\nu}_\ga  (\ft_\ga^{-1}(A) \cap E^{-k}).
\end{eqnarray}
By $\ft_\ga$ invariance of $\tilde{\nu}_\ga$,
$ \tilde{\nu}_\ga (\ft^{-1}_\ga(A) \cap  ([0,1]\backslash \Et_\ga)) $
equals
$ \tilde{\nu}_\ga (A \cap  ([0,1]\backslash \Et_\ga)) $
and
$\sum_{k=0}^l \tilde{\nu}_\ga  (\ft_\ga^{-1}(A) \cap E^{-k})$
equals
$\sum_{k=0}^l \tilde{\nu}_\ga  (A \cap E^{-k})$.
The right hand side of (\ref{eqn_inv}) therefore is $\bar{\nu}_\ga (A)$,
showing $f_\ga$-invariance of $\bar{\nu}_\ga$ and thus of $\nu_\ga$.
\qed \\

The following proposition serves several purposes, it is used 
to provide estimates for the average 
length of laminar and relaminarization phases.
Its proof will be postponed to Section~\ref{app}.

\begin{prop}\label{prop_outside}
Let $V$ be a small neighborhood of $c$.
For $x\in [0,1]$, let $\tilde{L}_{V}$ 
be the number of iterations under $\ft_\ga$
required for $x$ to enter $V$.
Then there is $L>0$ so that for all $\ga \in \Omega$,
\begin{eqnarray*}
\int_{[0,1]} \tilde{L}_{V} d\tilde{\nu}_\ga/ \tilde{\nu}_\ga ([0,1]\backslash V) 
& \le & L.
\end{eqnarray*}
\end{prop}

Corresponding statements for the study of the boundary crisis
bifurcation are contained in \cite{homyou00}.
Note that $V$ need not contain $\Et_\ga$; a similar statement
where $\ft_\ga$ is replaced by $f_\ga$ is therefore untrue.
However, Proposition~\ref{prop_outside} has the following corollary
which deals with $f_\ga$.
Paraphrasing, it shows that a 
typical (with respect to the invariant measure $\nu_\ga$) 
point $x \in [0,1]$ needs a bounded number of iterates to enter $\Et_\ga$.

\begin{prop}\label{prop_relaminar}
For $x\in [0,1]$, let $L_{\Et_\ga}$ be the number of iterations under $f_\ga$
required for $x$ to enter $\Et_\ga$.
Then there is $L>0$ so that for all $\ga \in \Omega$,
\begin{eqnarray*}
\int_{[0,1]} L_{\Et_\ga} d\nu_\ga/ \nu_\ga ([0,1]\backslash \Et_\ga) 
& \le & L.
\end{eqnarray*}
\end{prop}

\noindent {\sc Proof.}
At $\ga = 0$, $f_0$ is renormalizable: there is an interval $V$ 
containing $c$ in its interior and
$a$ in its boundary, so that $f^q_0(V) \subset V$.
By slightly extending $V$, we may assume 
that it contains $\Et_\ga$ 
for small values of $\ga$.
The result follows from Proposition~\ref{prop_outside}
by noting that outside $\Et_\ga$,
$\ft_\ga$ equals $f_\ga$. \qed \\

We will now show how Theorem~\ref{intermittency} is proved by
combining Propositions~\ref{prop_tildeacip},~\ref{prop_outside}
and~\ref{prop_relaminar}.\\

\noindent {\sc Proof of Theorem~\ref{intermittency}.}
Let $\Omega$ be as constructed in Section~\ref{sec_br}
and take $\ga \in \Omega$.
Let $\nu_\ga$ be the absolutely continuous invariant measure
for $f_\ga$, obtained in Lemma~\ref{prop_mb}.
The measure $\nu_\ga$ is ergodic, so that by Birkhoff's Ergodic Theorem,
\begin{eqnarray}\label{ergodic}
\lim_{m\to \infty} \sum_{i=0}^{m-1} 1_A (f^i_\ga(x)) & = & \nu_\ga (A),
\end{eqnarray}
for any Borel set $A \subset [0,1]$ and almost all $x\in [0,1]$.
It follows that also
\begin{eqnarray}\label{ergodictilde}
\lim_{m\to \infty} \sum_{i=0}^{m-1} 1_A (\ft^i_\ga(x)) & = & \tilde{\nu}_\ga (A).
\end{eqnarray}
Hence, for almost all $x \in [0,1]$, the distribution with which
points in the orbit $\{\ft_\ga^i(x) \}$ are in  $\Et_\ga$ is given by
$\tilde{\nu}_\ga$.
Let $I$ be a compact interval in $[d,a)$.
If $V$ is a neighborhood of $c$, then by (\ref{ergodictilde})
and applying Proposition~\ref{prop_outside}, the measure $\tilde{\nu}_\ga(V)$
of $V$ is bounded from below by a positive constant, uniformly in $\ga$. 
Observe that by invariance of $\tilde{\nu}_\ga$, 
$\tilde{\nu}_\ga ( \ft_\ga^k(V)) \ge \tilde{\nu}_\ga (V)$, $k \ge 0$.
Therefore, the measure $\tilde{\nu}_\ga(I)$
of $I$ is bounded from below by a positive constant, uniformly in $\ga$. 
An easy computation shows that  
for any compact interval $I$ inside $[d,a)$, 
the number of iterates needed for a point $x \in I$ to leave $\Et_\ga$
is bounded from below by $K /\sqrt{\ga}$ for some $K>0$. 
It follows that the average duration of orbit pieces
of $f_\ga$ in $\Et_\ga$ is bounded from below by
$K /\sqrt{\ga}$ for some $K>0$.
Combining this with Proposition~\ref{prop_relaminar}
proves the upper bound on $\chi_{\bar{E}}(\ga)$ in Theorem~\ref{intermittency}.
The lower bound is a trivial consequence of the fact that there is always
a positive number of iterates between two laminar phases and that the
maximum number of consecutive iterations in $\tilde{E}$ is bounded
above by $K'/\sqrt{\ga}$ for some constant $K'$.

The argument to show that $\nu_\ga$ converges weakly to $\nu_0$ is similar.
By definition of $\bar{\nu}_\ga$,
$$
\bar{\nu}_\ga (E^{-j}) = \sum_{k=j}^l \tilde{\nu}_\ga (E^{-k}).
$$
Reasoning as above one shows that $\tilde{\nu}_\ga(E^{-j})$, $0\le j \le l$,
is bounded from below, uniformly in $\ga$.
It follows that 
$\bar{\nu}_\ga (\Et_\ga)$ gets arbitrarily large as $\ga \to 0$
(because then $l\to \infty$).
Since $\bar{\nu}_\ga( [0,1]\backslash \Et_\ga) = 
\tilde{\nu}_\ga ( [0,1]\backslash \Et_\ga)$ is obviously bounded, 
$\nu_\ga (\Et_\ga) \to 1$ as $\ga \to 0$.
Because this holds for any neighborhood $\Et_\ga$, this shows that
$\nu_\ga \rightharpoonup \nu_0$ as $\ga \to 0$. 
\qed

\subsection{Relaminarization}\label{app}

Purpose of this section is to indicate a proof of Proposition~\ref{prop_outside}.
To introduce the reasoning, we start with an alternative proof for
Proposition~\ref{prop_relaminar}, which does not derive it as a corollary to
Proposition~\ref{prop_outside}.\\

\noindent {\sc Proof of Proposition~\ref{prop_relaminar}.}
We claim that $\tilde{\nu}_\ga (([0,1]\backslash \Et_\ga)$ 
is bounded away from 0, 
uniformly in $\ga$.
To establish the claim, note that $\ft^{-1} ([0,1]\backslash \Et_\ga)$ and
$[0,1]\backslash \Et_\ga$ together cover $[0,1]$. 
By $\ft_\ga$-invariance of $\tilde{\nu}_\ga$, the measures of both sets are the same,
and, since they add up to at least 1, are bounded from below by $\frac{1}{2}$.
The normalized measure $\nu_\ga / \nu_\ga ([0,1]\backslash \Et_\ga)$
equals $\tilde{\nu}_\ga / \tilde{\nu}_\ga ([0,1]\backslash \Et_\ga)$
outside $\Et_\ga$.
Therefore, 
applying Proposition~\ref{prop_tildeacip}, it follows that
a uniform bound of the form
\begin{eqnarray}\label{unifbound}
\nu_\ga (A) / \nu_\ga ([0,1]\backslash \Et_\ga) & \le & K \sqrt{ m(A)}
\end{eqnarray}
holds for
Borel sets $A \subset ([0,1]\backslash \Et_\ga)$.

\begin{figure}[ht]

\centerline{\epsfig{figure=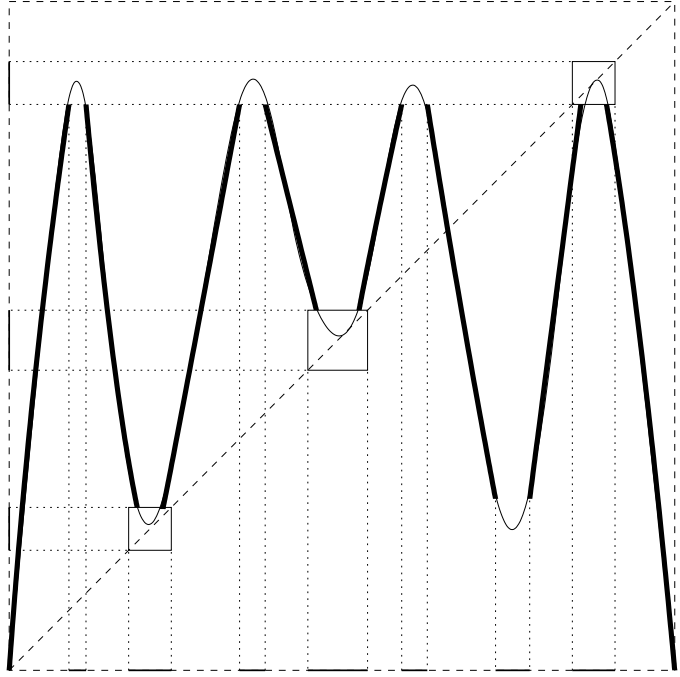,height=9cm,width=9cm}}

\caption{This figure illustrates the proof of Proposition~\ref{prop_relaminar}.
It shows the third iterate of a unimodal map at or near 
a saddle-node bifurcation.
Projected on the left side of the box are the intervals that make up $\bar{N}$.
Projected on the bottom one finds $O$. 
 \label{relaminar} }
\end{figure}

At $\ga = 0$, $f_0$ is renormalizable: there is an interval $V$ 
containing $c$ in its interior and
$a$ in its boundary, so that $f^q_0(V) \subset V$.
By slightly extending $V$, we may assume 
that it contains $\Et_\ga$ 
for small values of $\ga$.
Write $\bar{V} = V \cup f_\ga(V) \cup \cdots \cup f^{k-1}_\ga(V)$.
Write $O$ for the union of $\bar{V}$ and those intervals in 
$f^{-q}_\ga (\bar{V})$ 
that contain a critical
point for $f^q_\ga$, see Figure~\ref{relaminar}.
Because any point in $O \setminus \Et_\ga$ either enters $\Et_\ga$ 
in a uniformly bounded number of iterates ($q$),
it suffices to prove the proposition with $\Et_\ga$ replaced by $O$.

Of course, also for Borel sets $A \subset [0,1]\backslash O$ a uniform bound of the form 
$\nu_\ga (A) / \nu_\ga ([0,1]\backslash O) \le K \sqrt{ m(A)}$ holds.
Let $h_\ga: [0,1]\backslash O \to O$ denote the restriction of
$f_\ga$ to $[0,1]\backslash O$.
By \cite{piayor79}, $h_\ga$ possesses a conditionally invariant measure $\zeta_\ga$;
$\zeta_\ga$ is characterized by
$$
\zeta_\ga(A) = \zeta_\ga (h^{-1}_\ga (A)) / \zeta_\ga (h^{-1}_\ga ([0,1]\backslash O)
$$
for Borel sets $A$.
Moreover, $\zeta_\ga$ is absolutely continuous
with respect to Lebesgue measure and has density bounded 
and bounded away from 0 on its support, 
uniformly in $\ga$.
Say
$$ 
\frac{1}{C} m(A) \le \zeta_\ga (A)
$$
for Borel sets $A$.
Write $(h_\ga)_* \zeta_\ga = \beta_\ga \zeta_\ga$ for some $0<\beta_\ga<1$.
Since $m(O)$ is bounded from below, $\beta_\ga$ is bounded away from 1.
Applying (\ref{unifbound}),
\begin{eqnarray*}
\nu_\ga (h_\ga^{-n}([0,1])) & \le & K \sqrt{m(h_\ga^{-n}([0,1]))}
\\
& \le & K \sqrt{C} \sqrt{\zeta_\ga(h_\ga^{-n} ([0,1]))}
\\
& = & 
K \sqrt{C} \sqrt{\beta_\ga}^n.
\end{eqnarray*}
It follows that
$$
\int_{[0,1]} L_O d\nu_\ga = \sum_{i\ge 1} \nu_\ga (h_\ga^{-i}([0,1]))
$$
is bounded uniformly in $\ga$.
\qed \\

\noindent {\sc Proof of Proposition~\ref{prop_outside}.}
Let $N$ be the interval containing $c$ so that $f_0^q (N) = N$.
We may assume that $V \subset N$ and $V \cap \Et_\ga = \emptyset$.
Given $V$, let $\bar{V} = \cup_{i=0}^{q-1} f^{-i}_\ga(V) \cap f^{q-i}_\ga (N)$
be the union of the interval in $\bar{N}$ that are mapped onto $V$ by 
iterates of $f_\ga$. 
Write $O$ for the union of $\bar{V}$ and those intervals in 
$f^{-q}_\ga (\bar{V})$ 
that contain a critical
point for $f^q_\ga$. 
Restricted to $\bar{V}$, $f_\ga$ equals $\ft_\ga$.
Include into $O$ the set of points $\ft^{-q}_\ga (\ft_\ga (D))$, where
$D$ is the set of discontinuities of $\ft_\ga$.
Write $\tilde{h}_\ga$ for the restriction of $\ft_\ga$ to $[0,1]\backslash O$.
By Proposition~\ref{prop_exp}, some iterate of $\tilde{h}_\ga$ 
is expanding.
Observe that the number of branches of $\tilde{h}_\ga$ is constant
for $\ga \in [\ga_{l+1},\ga_l]$, but increases with $l$ since
the branches of $\Res{ \ft_\ga}{\Et_\ga}$ are included.  
For fixed $l$, 
Proposition~\ref{prop_outside} follows for $\ga \in [\ga_{l+1},\ga_l] \cap \Omega$
as in the above indicated alternative proof of Proposition~\ref{relaminar}.
To obtain uniform bounds in $l$, we must investigate properties of
conditionally invariant measures for $\tilde{h}_\ga$.
Following \cite{piayor79}, there is a conditionally invariant measure
$\tilde{\zeta}_\ga$ for $\tilde{h}_\ga$;
$\tilde{\zeta}_\ga$ is characterized by 
$$
\tilde{\zeta}_\ga (A) = \tilde{\zeta}_\ga (\tilde{h}_\ga^{-1}(A)) / \tilde{\zeta}_\ga(\tilde{h}_\ga^{-1}([0,1]\backslash O).
$$
The proof of Proposition~\ref{prop_outside} is identical to the above proof of 
Proposition~\ref{prop_relaminar}, once Lemma~\ref{lemma} is proved. \qed \\

\begin{lem}\label{lemma}
Let $\tilde{h}_\ga$ be as in the above proof of Proposition~\ref{prop_outside}.
The density of the conditionally invariant measure $\tilde{\zeta}_\ga$
for $\tilde{h}_\ga$ is bounded away from 0
uniformly in $\ga$.
\end{lem}

\noindent {\sc Proof.}
By Proposition~\ref{prop_exp}, an iterate $\tilde{h}_\ga^N$ of $\tilde{h}_\ga$ is an expansion.
Hence, the existence \cite{piayor79} of the conditionally invariant measure follows
from \cite{piayor79}. To derive bounds on its density, we must examine 
the existence proof. 

The conditionally invariant measure $\tilde{\zeta}_\ga$ is constructed by finding
its density as a fixed point of a Perron-Frobenius operator.
Write ${\cal C}([0,1])$ for the set of positive continuous functions $g$ on $[0,1]$
with $\int_{[0,1]} g dm =1$.
Near a point $x \in [0,1]$, $\tilde{h}_\ga$ has a number of 
inverse functions $\psi_i$.
Define a Perron-Frobenius operator $P_\ga$ on ${\cal C}([0,1])$ by
$$
P_\ga g (x) = \sum_i | \psi_i'(x)| g \circ \psi_i (x) / \int_{\tilde{h}^{-1}_\ga ([0,1])} g dm
$$
A fixed point of $P_\ga$ is the density of a conditionally invariant measure.
It is shown in \cite{piayor79} that 
$P_\ga$ possesses a unique fixed point $\xi_\ga$ and that for each
$g \in {\cal C}([0,1])$,
\begin{equation}\label{converges}
\lim_{n\to\infty} P_\ga^n g  = \xi_\ga.
\end{equation}
The fixed point $\xi_\ga$ is the density of the measure $\tilde{\zeta}_\ga$.

For a Lipschitz density $g$, let the regularity of $g$ be given by
$$
\mbox{Reg } g = \sup \{|g'(x)|/g(x); \;\; x\in [0,1], g'(x)\mbox{ is defined}, g(x)>0 \}.
$$
We will show  that
\begin{equation}\label{boundedregularity}
\limsup_{n\to \infty} \mbox{Reg } P_\ga^n g \le \rho,
\end{equation}
for some $\rho$ independent of $g$.
We remark that from \cite{piayor79}
one concludes that such a bound holds when $\ga$ is restricted to 
an interval $[\ga_{l+1},\ga_l]$ (see also \cite{lasyor81}).
Their arguments do not imply (\ref{boundedregularity}),
the difficulty being the growth of the number of branches of $\tilde{h}_\ga$ as $\ga\to 0$.

To evaluate $\mbox{Reg }P_\ga g$, compute
\begin{eqnarray*}
\frac{ | (P_\ga g)'(x)|}{P_\ga g(x)} & = & \frac{ \left| \sum_i  [ \psi_i'(x) g\circ \psi_i (x)]' \right| }
                                      { \sum_i  \psi_i'(x) g\circ \psi_i (x)}
\\
& \le &
\frac{ \left| \sum_i   \psi_i''(x) g\circ \psi_i (x) \right|}
     { \sum_i  \psi_i'(x) g\circ \psi_i (x)}
+
\frac{ \left| \sum_i  \psi_i'(x) g'\circ \psi_i (x) \psi_i'(x) \right|}
     { \sum_i  \psi_i'(x) g\circ \psi_i (x)}
\\
& \le &
\max_i \frac{| \psi_i''(x)|}{\psi_i'(x)} + 
\max_i |\phi_i'(x)| \frac{g'\circ \psi_i (x)}{g\circ \psi_i(x)}.
\end{eqnarray*}
In the last step, we used that 
$| \sum a_i / \sum b_i | \le \max_i | a_i/b_i |$ for numbers $a_i$ and positive numbers $b_i$.
We claim that $\max_i \frac{| \psi_i''(x)|}{\psi_i'(x)}$ is bounded, uniformly in $\ga$.
This is clear for $x$ outside $\Et_\ga$, where $\psi_i'$ is bounded away from 0 and
$\psi_i''$ is bounded.
Take $x \in \Et_\ga$ and consider the inverse branch $\psi_i$ provided by
the flow $\phi_\ga^t$ of the adapted vector field $\phi_\ga$ (see Proposition~\ref{embedding}).
Thus $ \frac{d}{dt} \phi_\ga^t (x_0) = \phi_\ga ( \phi_\ga^t(x_0) )$.
To bound $|\psi_i''| / \psi_i'$ we consider the flow
$\phi_\ga^t$ for negative $t$. 
Let $J_\ga^t = (\phi_\ga^t)'' / (\phi_\ga^t)'$.
Then $J_\ga^t$ is a solution of 
\begin{eqnarray}\label{J}
\frac{d}{dt} J_\ga^t (x_0) = \Res{\phi_\ga''}{\phi_\ga^t (x_0)} (\phi_\ga^t)' (x_0).
\end{eqnarray}
For $\ga$ small, $\phi_\ga''$ is close to $\phi_0''$, which is nonzero by assumption.
Further, $(\phi_\ga^t)' (x_0)$ is a solution of 
$\frac{d}{dt} (\phi_\ga^t)' (x_0) = 
\Res{ \phi_\ga'}{\phi_\ga^t(x_0)}  (\phi_\ga^t)' (x_0)$.
Using this one can show that  $J_\ga^t$ is bounded, uniformly in $\ga$.
We sketch a possible reasoning.
One may assume that $\phi_\ga (x) = x^2 + \ga + {\cal O}(x^2 \ga) + {\cal O}(x^3)$.
Let $\phi_{\ga,a}(x)  = a x^2 + \ga$ and denote the corresponding flow by
 $\phi_{\ga,a}^t$. 
For arguments sake, take $x_0 < 0$.
Then, for $t < 0$,  
$\phi_{\ga,a^+}^t (x_0) \le \phi_\ga^t (x_0) \le \phi_{\ga,a^-}^t (x_0)$ 
for some $a_- < a_+$ close to 1.
Solving for  $\phi_{\ga,a}^t$ shows that
\begin{eqnarray}\label{diffeqn}
\frac{1}{\sqrt{a}} \arctan (\sqrt{\frac{a}{\ga}} \phi_{\ga,a}^t (x_0)) & = &
\sqrt{\ga} t + \frac{1}{\sqrt{a}} \arctan (\sqrt{\frac{a}{\ga}} x_0).
\end{eqnarray}
Note that  
$\frac{d}{dt}  (\phi_{\ga,a}^t)'  (x_0) = 2 a \phi_{\ga,a}^t  (x_0)  (\phi_{\ga,a}^t)'  (x_0)$.
It follows that
$$
 (\phi_{\ga,a^-}^t)'  (x_0)\le
 (\phi_{\ga}^t)' (x_0) \le 
 (\phi_{\ga,a^+}^t)'  (x_0).
$$
Differentiating (\ref{diffeqn}) yields
$$
(\phi_{\ga,a}^t)'  (x_0) =  \frac{ a (\phi_{\ga,a}^t (x_0))^2 + \ga}{ a x_0^2 + \ga}.
$$
These bounds and (\ref{J}) prove that $J_\ga^t$ is bounded, uniformly in $\ga$.

Replacing $\tilde{h}_\ga$ by $\tilde{h}_\ga^N$, which is an expanding map, 
the above reasoning shows that
$$
\mbox{Reg }P_\ga^N g \le M + \frac{1}{\lambda} \mbox{Reg }g
$$
for constants $M>0$ and $\lambda > 1$.
Now (\ref{boundedregularity}) follows by iterating the bound on the regularity of $P_\ga^N g$.

We can now basically follow the proof of Proposition~7 in  \cite{piayor79} to
conclude the lemma.
Take $g \in {\cal C}([0,1])$ with $\mbox{Reg }g \le \rho$.
Then, for $x \in [0,1]$, 
$$
P_\ga^n g (x) = \sum_i | \psi_i'' (x)|  g\circ \psi_i (x),
$$
where $\psi_i$ are the inverse branches of $\tilde{h}_\ga^n$ near $x$. 
Since $\int_{[0,1]} g dm = 1$, there is  
a component $A_j$ of $[0,1]\backslash O$ on which $\sup_{A_j} g \ge 1$.
We may also assume that $\sup_{A_j} g \ge \frac{1}{2}$ on a component $A_j$
with $m(A_j) \ge \beta$ for some constant $\beta > 0$.
Then $\inf_{A_j}g \ge \frac{1}{2} e^{-\beta \rho}$.

If $n$ is large enough, $[f_\ga^2(c),f_\ga(c)] \subset \tilde{h}_\ga^n (A_j)$ for all
small $\ga$.
Hence, for each $x \in [f_\ga^2(c),f_\ga(c)]$, there exists $i_0$ with $\psi_{i_0} (x) \in A_j$.
This gives
\begin{eqnarray*}
P^n_\ga g (x) & \ge & |\psi_{i_0}'' (x)| g\circ \psi_{i_0} (x),
\\
& \ge &
|\psi_{i_0}'' (x)| e^{-\rho \beta}. 
\end{eqnarray*}
Since $m(A_j) \ge \beta$, we may take $i_0$ so that $|\psi_{i_0}'' (x)|$ is bounded from below.
Therefore,
$P^n_\ga g (x) \ge d$ for some $d>0$ which is independent of $g$ and $\ga$.
By (\ref{converges}), this proves the lemma.
\qed\\

\end{document}